\documentclass{amsart}
\usepackage{amssymb}

\nonstopmode

\numberwithin{equation}{section}
\newtheorem{theorem}{Theorem}[section]
\newtheorem{lemma}[theorem]{Lemma}
\newtheorem{proposition}[theorem]{Proposition}
\newtheorem{corollary}[theorem]{Corollary}

\theoremstyle{definition}
\newtheorem{definition}[theorem]{Definition}
\newtheorem{remark}[theorem]{Remark}
\newtheorem{example}[theorem]{Example}

\newcommand\codim{\operatorname{codim}}
\newcommand\coker{\operatorname{coker}}
\newcommand\Ann{\operatorname{Ann}}
\newcommand\Supp{\operatorname{Supp}}

\newcommand\Hom{\operatorname{Hom}}
\newcommand\Ext{\operatorname{Ext}}

\newcommand\rank{\operatorname{rank}}
\newcommand\rankk{\operatorname{rank}_K}
\newcommand\depth{\operatorname{depth}}

\newcommand\im{\operatorname{im}}

\newcommand{\HH}{H_{\mathfrak m}}

\newcommand{\TR}{\operatorname{Tor}^R}
\newcommand{\ER}{\operatorname{Ext}_R}
\newcommand{\Proj}{\operatorname{Proj}}
\newcommand{\Rad}{\operatorname{Rad}}

\newcommand{\mif}{\mbox{if} ~}

\newcommand{\ffi}{\varphi}

\newcommand{\cJ}{{\mathcal J}}
\newcommand{\cB}{{\mathcal B}}
\newcommand{\cC}{{\mathcal C}}
\newcommand{\cD}{{\mathcal D}}
\newcommand{\cE}{{\mathcal E}}
\newcommand{\cF}{{\mathcal F}}
\newcommand{\cG}{{\mathcal G}}
\newcommand{\cO}{{\mathcal O}}
\newcommand{\cP}{{\mathcal P}}
\newcommand{\cQ}{{\mathcal Q}}

\newcommand{\cN}{{\mathcal N}}

\newcommand{\cBf}{{\mathcal B}_{\varphi}}

\newcommand{\Bf}{B_{\varphi}}
\newcommand{\Mf}{M_{\varphi}}

\newcommand{\fm}{{\mathfrak m}}

\newcommand{\cOP}{\Omega_{\mathbb{P}^n}} 
\newcommand{\cPP}{{\mathcal O}_{\mathbb{P}^n}}

\newcommand {\ZZ}{\mathbb{Z}}

\newcommand {\PP}{\mathbb{P}}

\newcommand {\ds}{\subset\! \! \! \! + \,}

\begin{document}
\title[Buchsbaum-Rim sheaves]{Buchsbaum-Rim sheaves and their multiple
  sections} 
\author[Juan C.\ Migliore, Uwe Nagel and Chris Peterson]{Juan C.\
  Migliore$^1$, Uwe Nagel$^2$ and Chris Peterson$^3$} 

 
\thanks{$^1$Department of Mathematics, 
        University of Notre Dame, 
        Notre Dame, IN 46556, 
        USA, Juan.C.Migliore.1@nd.edu \\
$^2$Fachbereich Mathematik und Informatik,
Universit\"at-Gesamthochschule Paderborn,  D--33095 Paderborn, Germany,
uwen@uni-paderborn.de \\
$^3$Department of Mathematics, 
        University of Notre Dame, 
        Notre Dame, IN 46556, 
        USA, peterson@math.nd.edu}

\begin{abstract} This paper  begins by introducing  and characterizing
  Buchsbaum-Rim sheaves 
  on $Z = \Proj R$ where $R$ is a graded Gorenstein $K$-algebra. They are
  reflexive sheaves arising as the sheafification of kernels of sufficiently
  general maps between free $R$-modules. 

Then we study multiple sections of a Buchsbaum-Rim sheaf $\cBf$, i.e, we
consider morphisms $\psi: \cP \to \cBf$ of sheaves on $Z$
dropping rank in the expected codimension, where $H^0_*(Z,\cP)$ is a free
$R$-module. The main purpose of this paper
is to study  properties of schemes associated to the degeneracy locus $S$
of $\psi$. It turns out that $S$ is often not equidimensional. Let $X$
denote the top-dimensional part of $S$. In this paper we measure the
``difference'' between $X$ and $S$, compute their cohomology modules and
describe ring-theoretic properties of their coordinate rings. Moreover, we
produce graded free resolutions of $X$ (and $S$) which are in general
minimal. 

Among the applications we show how one can embed a subscheme into an
arithmetically Gorenstein subscheme of the same dimension and prove that
zero-loci of sections of the dual of a null correlation bundle are
arithmetically Buchsbaum. 
\end{abstract}

\maketitle 

\tableofcontents 

\section{Introduction}

A fundamental method for constructing algebraic varieties is to
consider the degeneracy locus of a morphism between a pair of coherent
sheaves. By varying
the morphism one obtains families of varieties. By placing various
restrictions on the
coherent sheaves one can force the degeneracy locus to have special
properties. These extra
restrictions may also provide for more tools with which to study the
degeneracy locus. If many
restrictions are placed on the sheaves then a great deal of precise
information can be
extracted but at the expense of generality. If one puts no restrictions on
either of the
sheaves then, of course, very little information can be extracted
concerning the degeneracy
locus. In this paper we take the middle road between these two extremes. We
consider a class
of sheaves, which while quite general, behave well enough that a substantial
amount of
information can be obtained with respect to their degeneracy loci. The only
restriction placed
on the sheaves is that they arise as the sheafification of the kernels of
sufficiently
general maps between free $R$-modules, where $R$ is a Gorenstein
$K$-algebra.  We will refer to
the sheaves constructed in this manner as Buchsbaum-Rim sheaves.

The purpose of
this paper is to introduce the class of Buchsbaum-Rim sheaves, to elicit their
main properties,
and to make a systematic and detailed study of the degeneracy loci obtained
by taking
multiple sections of these sheaves. Although Buchsbaum-Rim sheaves are
necessarily reflexive (they are the sheafification of a second syzygy
module) they will
not, in general, be locally free. A number of tools are used to manipulate
and control these
objects. Certainly techniques developed by Eagon-Northcott, Buchsbaum-Rim,
Buchsbaum-Eisenbud,
Kirby and Kempf play an important role. These are combined with the methods
of local cohomology
and several homological techniques to produce the main results.
The paper opens with a brief section providing preliminary background
information of use in
later sections of the paper. Sections three, four and five form the
technical heart and the
final section closes out the paper with three applications.

Section three introduces
Buchsbaum-Rim sheaves and Buchsbaum-Rim modules. To begin with we should
make clear the
definition of a Buchsbaum-Rim sheaf. In the following $R$ will always
denote a graded Gorenstein $K$-algebra of dimension $n+1$ where $K$ is an
infinite field. Furthermore, the scheme $Z$ will be the projective spectrum
of $R$. 

\begin{definition} Let $\cF$ and $\cG$ be locally free sheaves of ranks $f$
and $g$
respectively on $Z$.  Let $\ffi: \cF \to \cG$ be a generically surjective
morphism. Suppose that the degeneracy locus of
  $\ffi$ has codimension $f-g+1$ and that the modules $F =
  H^0_*(Z,\cF)$ and $G = H^0_*(Z,\cG)$ are  free
  $R$-modules. We call the kernel of $\ffi$ a {\it Buchsbaum-Rim sheaf}
and denote it by $\cBf$.  

By abuse of notation we denote the homomorphism $F \to G$ induced by $\ffi$
again by $\ffi$. Moreover, we put $\Bf = H^0_*(Z,\cBf)$ and $\Mf = \coker
\ffi$, so that we have an exact sequence
$$
0 \to B_{\ffi} \to F \stackrel{\ffi}{\longrightarrow} G \to  M_{\ffi} \to 0.
$$
We call $\Bf$ a {\it Buchsbaum-Rim module}.
\end{definition} 

Thus the cotangent bundle of projective space is a Buchsbaum-Rim sheaf. 
The cohomology of its exterior powers is given by the Bott formula. Letting
$r(R)$ denote the index
of regularity of a graded ring we have the following lemma which
generalizes the Bott formula
to arbitrary Buchsbaum-Rim sheaves. 

\begin{lemma} Let
  $\Bf$ be a Buchsbaum-Rim module of rank $r$. Then it holds:
\begin{itemize}

\item[(a)] For $i = 0,\ldots r$ there are isomorphisms
$$
\wedge^{r-i} \cBf^* \otimes \wedge^f \cF \otimes \wedge^g \cG^* \cong (\wedge^i
\cBf^*)^*. 
$$
\item[(b)] For $i < r$ we have
$$
H^j_*(Z, \wedge^i \cBf^*) \cong \left \{ \begin{array}{ll}
0 & \mif 1 \leq j \leq n \; \mbox{and} \; j \neq n-i \\
S_{i}(\Mf)^{\vee}(1 - r(R)) & \mif j = n-i.
\end{array} \right.
$$
\end{itemize}
\end{lemma}

The value of this lemma will become clear when we construct the generalized
Koszul complexes
associated to taking multiple sections of Buchsbaum-Rim sheaves. The lemma
also suggests a relationship to Eilenberg-MacLane sheaves. Recall that an
$R$-module $E$ is
called an {\it
  Eilenberg-MacLane module} of depth~$t$, $0 \leq t \leq n+1$ if
$\HH^j(E) = 0 \quad \mbox{for all} \;  j \neq t \; \mbox{where} \; 0 \leq j
\leq n.$
Similarly, a sheaf $\cE$ on $Z$ is said to be an {\it Eilenberg-MacLane
  sheaf} if $H^0_*(\cE)$ is an Eilenberg-MacLane module. A cohomological
characterization of
Buchsbaum-Rim sheaves can then be stated as follows.

 \begin{proposition}  A sheaf $\cE$ on $Z$
  is a Buchsbaum-Rim sheaf if and only if $E = H^0_*(\cE)$ is a reflexive
  Eilenberg-MacLane module with finite projective dimension and rank $r
  \leq n$ such that $\HH^{n-r+1}(E)^{\vee}$ is a perfect $R$-module of
  dimension  $n-r$ if $r \geq 2$.
\end{proposition}

In section four we make a detailed study of the cohomology of the
degeneracy locus of
multiple sections of Buchsbaum-Rim sheaves. Consider a morphism $\psi: \cP
\to \cBf$ of
sheaves of rank $t$ and $r$ respectively on  $Z = Proj (R)$,
where $\cBf$ is a Buchsbaum-Rim sheaf and $H^0_*(Z, \cP)$ is a free
$R$-module. If $t = 1$ then $\psi$ is just a section of some twist of
$\cBf$. For arbitrary $t<r$ we say that $\psi$ corresponds to taking
multiple sections of
$\cBf$. We always suppose that the degeneracy locus $S$ of $\psi$
has (the expected) codimension $r - t + 1 \geq 2$ in $Z$ (if $t=1$ then $S$
is just the zero-locus of a regular section of a Buchsbaum-Rim sheaf). 

 An Eagon-Northcott complex involving $\cBf$ will play an essential
role. Our approach will be algebraic and uses local cohomology. Taking
global sections we obtain an $R$-homomorphism
$$
\psi : P \to \Bf
$$
where $P$ is a free $R$-module of rank $t, 1 \leq t < r$. Then there is an
Eagon-Northcott complex 
$$
E_{\bullet}  \colon \quad 0 \to E_r \stackrel{\delta_r}{\longrightarrow} E_{r-1}
\stackrel{\delta_{r-1}}{\longrightarrow} \cdots \to E_t
\stackrel{\delta_t}{\longrightarrow} I(\psi) \otimes \wedge^t P^*  \to 0
$$
where
$$
E_i = \wedge^i \Bf^* \otimes S_{i - t}(P) 
$$
and $I(\psi)$ is the ideal defined by the image of $\delta_t = \wedge^t
\psi^*$. The saturation of $I(\psi)$ is the homogeneous ideal of the
degeneracy locus $S$.

With the help of this Eagon-Northcott complex and the first lemma we  get
the following formula 
for the cohomology modules of $S$. 

\begin{proposition} Let $I = I(\psi)$.  Then it holds for   $j\neq \dim R/I
  = n+t-r$ 
$$
H^j_m (R/I) \cong \left\{ \begin{array}{ll}
S_i(M_{\ffi})^\vee \otimes
S_{i-t} (P) \otimes \wedge^t P (1-r(R)) & \mbox{if } j=n+t-2i \quad
\mbox{where }
\max \{t, \frac{r+1}{2}\} \leq i \leq \left \lfloor \frac{r+t}{2} \right
\rfloor  \\
0 & \mbox{otherwise}. 
\end{array} \right.
$$
\end{proposition} 

This proposition allows us to decide if $S$ is equidimensional. It will
often turn out that this is not the case.  Thus we are
also interested in the top-dimensional component $X$ of $S$, i.e. the union
of the
highest-dimensional components of $S$. Let $J = J(\psi)$ denote the
homogeneous ideal of $X$. 
We need a measure of
the failure of $I$ to be equidimensional and we need a close relation
between the cohomology
of the schemes $X$ and $S$. This is provided in the following result.

\begin{proposition}  Letting $I$ and $J$ denote the ideals associated to
  $\psi$  we have:
 \begin{itemize}
\item[(a)]  $I$ is unmixed if and only if $r+t$ is odd.
\item[(b)] If $r+t$ is even then $I$ has a primary component of codimension
  $r+1$. Let $Q$ be the intersection of all such components. Then we have $I =
  J \cap Q$ and
$$
H^j_m(R/J) \cong \left\{ \begin{array}{ll}
H^j_m (R/I) & \mbox{if } j\neq n-r \\
0 & \mbox{if } j=n-r
\end{array}\right.
$$
and
$$
J/I \cong S_{\frac{r-t}{2}} (M_{\ffi}) \otimes S_{\frac{r-t}{2}}(P) \otimes
\wedge^f F^* \otimes \wedge^g G \otimes
\wedge^t P.
$$
\end{itemize}
\end{proposition}

Combining these propositions in the proper manner we are in a position to
prove the main
theorem of section four.

\begin{theorem}
With the  notation above we have:
\begin{itemize}
\item[(a)] If $r = n$ then $S$ is equidimensional and locally
  Cohen-Macaulay.
\item[(b)] $S$ is equidimensional if and only if $r+t$
  is odd or $r = n$.
 Moreover, if $r < n$ then $X$ is locally Cohen-Macaulay if and only if $X$ is
  arithmetically Cohen-Macaulay.
\item[(c)] If $r+t$ is odd then $X = S$ is arithmetically Cohen-Macaulay if
  and only if $t=1$. In this case $S$ has Cohen-Macaulay type $\leq 1 +
  \binom{\frac{r}{2}+g-1}{g-1}$.
\item[(d)] Let  $r+t$ be  even. Then
\begin{itemize}
\item[(i)] $X$ is arithmetically Cohen-Macaulay if and only if $1 \leq t
\leq 2$. If $t =1$ then $X$ is arithmetically Gorenstein. If $t = 2$ then $X$
has Cohen-Macaulay type $\leq r - 1 + \binom{\frac{r}{2}+g-1}{g-1} \cdot
(\frac{r}{2} - 1)$.
\item[(ii)] If in addition  $r < n$ then the components
  of $S$ have either codimension $r-t+1$ or codimension $r+1$.
\end{itemize}
\end{itemize}
\end{theorem}

Contained in the theorem is the following surprising conclusion. Let $\cBf$
be an odd rank
Buchsbaum-Rim sheaf and let $X$ denote the top dimensional component of the
zero-locus of
any regular section of $\cBf$. Then $X$ is arithmetically Gorenstein. This
generalizes the main theorem of \cite{Mig-P_gorenstein}, where the case
$r=3$ was considered. 
\smallskip

Section five treats the problem of finding free resolutions of the
degeneracy loci. In order
to do this it is important to understand the homology modules of the
Eagon-Northcott complex $E_{\bullet}$ 
associated to $\psi$. We show that the homology can be summarized as follows.

\begin{proposition} The homology modules of the
  Eagon-Northcott $E_{\bullet}$ complex are:
$$
H_i (E_{\bullet}) \cong \left\{ \begin{array}{ll}
S_j(\Mf) \otimes S_{r-t-j}(P) \otimes \wedge^f F^* \otimes \wedge^g G &
\mif  i=r-1-2j
\mbox{ where } j \in \ZZ, \; t \leq i \leq r-3 \\
0 & \mbox{otherwise}.
\end{array}\right.
$$
\end{proposition}

This result allows us to conclude that $X$ and $S$ have  free resolutions
of finite length. However, it does not, in
general, provide
enough information to compute a minimal free resolution. To do this we need
several ingredients. First we need to understand the cohomology of the dual
of the Eagon-Northcott complex. This needs to be mixed with knowledge of
how the canonical
modules of $S$ and $X$ relate. The cohomology of the dual of the
Eagon-Northcott complex
is summarized in the following lemma, where $K_{R/I}$ denotes the canonical
module of $R/I(\psi)$. 

\begin{lemma} The dual of the Eagon-Northcott
  complex $E_{\bullet}$ provides a complex
$$
0 \to \wedge^t P \to E_t^* \stackrel{\delta^*_{t+1}}{\longrightarrow}
\ldots \to
  E_{r-1}^* \stackrel{\delta^*_{r}}{\longrightarrow} E_r^*
  \stackrel{\gamma}{\longrightarrow} K_{R/I} \otimes \wedge^t P (1 - r(R))
  \to 0
$$
which we denote (by slight abuse of notation) by $E_{\bullet}^*$. Its
(co)homology modules are given by
$$
H^i(E_{\bullet}^*) \cong \left \{
\begin{array}{ll}
S_j(\Mf) \otimes S_{j-t}(P)^* & \mif 2t + 1 \leq i = 2 j + 1 \leq r+1 \\
0 & \mbox{otherwise}
\end{array}
\right.
$$
In particular, $E_{\bullet}^*$ is exact if $t \geq \frac{r+1}{2}$.
\end{lemma}

To utilize these results on the Eagon-Northcott complex one still needs to
know when certain
terms in a free resolution can be split off. To do this we prove the
following result which, although rather
technical in appearance,  provides a substantial generalization to an
often-used result of
Rao. It can be applied in situations far removed from those addressed in
this paper and can
even be utilized when the ring $R$ is not a Gorenstein ring but is only a
Cohen-Macaulay ring.

\begin{proposition} Let $N$ be a finitely
  generated graded torsion $R$-module which has projective
  dimension~$s$. Then it holds for all integers $j \geq 0$  that 
  $\TR_{s-j}(N,K)^{\vee}$ is a direct summand of
$$
\oplus_{i=0}^j \TR_{j-i}(\ER^{s-i}(N,R),K).
$$
Moreover, we have $\TR_{s}(N,K)^{\vee} \cong \TR_{0}(\ER^{s}(N,R),K)$ and
that $\TR_{1}(\ER^{s}(N,R),K)$ is a direct summand of
$\TR_{s-1}(N,K)^{\vee}$.
\end{proposition}

Together these results allow us to write down a free resolution of the
degeneracy locus which is in general minimal. 
Thus the main theorem of section five is the
following, which gives the
free resolution for the degeneracy locus of a morphism $\psi: \cP \to
\cBf$, where
$\cP$ has rank $t$ and $\cBf$ has rank $r$.

\begin{theorem}  Consider the following
  modules where we use the conventions that $i$ and $j$ are non-negative
  integers and that a sum is trivial if it has no summand:
$$
A_k = \bigoplus_{\begin{array}{c}
{\scriptstyle i+2j = k + t -1}\\ [-4pt]
{\scriptstyle t \leq i+j \leq \frac{r+t-1}{2}}
\end{array}}
\wedge^i F^* \otimes S_j(G)^* \otimes S_{i+j-t}(P),
$$
$$
C_k = \bigoplus_{\begin{array}{c}
{\scriptstyle i+2j = r+1-t-k}\\ [-4pt]
{\scriptstyle i+j \leq \frac{r-t}{2}}
\end{array}}
\wedge^i F \otimes S_j(G) \otimes S_{r-t-i-j}(P) \otimes \wedge^f F^*
\otimes \wedge^g G.
$$
Observe that it holds:

$A_r = 0$ if and only if $r+t$ is even,

$C_1 = 0$ if and only if $r+t$ is odd,

$C_k = 0$ if $k \geq r+2-t$  and

$C_{r+1-t} = S_{r-t}(P) \otimes \wedge^f F^* \otimes \wedge^g G$. \\
Then the homogeneous ideal $I_X = J(\psi)$ of the top-dimensional part $X$
of the degeneracy locus $S$ has a graded free resolution of the form
$$
0 \to A_r \oplus C_r \to \ldots \to A_1 \oplus C_1 \to I_X \otimes \wedge^t
P^* \to 0.
$$
\end{theorem}

 Note that previously minimal free resolutions were
known only for a few classes besides the determinantal ideals.  
A number of examples of particular interest round out the section and
illustrate the theorem. 
\smallskip

The final section gives several additional applications that may be of
independent interest.
We show how Buchsbaum-Rim sheaves can be used to situate arbitrary
equidimensional schemes of
arbitrary codimension into arithmetically Gorenstein schemes. This will be
of relevance when one considers the
problem of linkage by arithmetically Gorenstein ideals as opposed to
complete intersection
linkage theory. We also show how to utilize Buchsbaum-Rim sheaves to
produce interesting new
examples of $k$-Buchsbaum sheaves as well as of arithmetically Buchsbaum
schemes. Finally we construct new vector bundles of rank $n-1$ on $\PP^n$
if $n$ is odd. We call them generalized null correlation bundles and show
that our results apply to  multiple sections of their duals.


\section{Preliminaries} \label{preliminaries}

Let $R$ be a  ring. If $R = \oplus_{i \in \mathbb{N}} R_i$ is graded then the
irrelevant maximal ideal $ \oplus_{i > 0} R_i$ of $R$ is denoted by
${\mathfrak m}_R$ or simply ${\mathfrak m}$. It is always assumed that
$R_0$ is an infinite field $K$ and that the $K$-algebra $R$ is generated by
the elements
of $R_1$. Hence $(R,\fm)$ is $^*$local in the sense of \cite{Bruns-Herzog}.

If $M$ is a module over the graded ring $R$ it is always assumed to be
$\mathbb{Z}$-graded. The set of its homogeneous elements of degree $i$ is
denoted by $M_i$ or $[M]_i$. All homomorphisms between graded $R$-modules
will be morphisms in the category of graded $R$-modules, i.e., will be
graded of degree zero. If $M$ is assumed to be a graded $R$-module it is
always understood that $R$ is a graded $K$-algebra as above.  We refer to the
context just described as the graded situation.
Although we are mainly interested in graded objects we note that our
results hold also true (with the usual modifications) in a local
situation. Then $(R,{\mathfrak m})$ will
denote a local ring with maximal ideal ${\mathfrak m}$.

If $M$ is an $R$-module, $\dim M$ denotes the Krull dimension of $M$. The
symbols $\rank_R$ or simply $\rank$ are reserved to denote the rank of $M$
in case it has one. For a $K$-module, $\rankk$ just denotes
 the vector space dimension over the field $K$.

There are two types of duals of an $R$-module $M$ we are going to use. The
$R$-dual of $M$ is $M^* =  \Hom_R(M,R)$. If $M$ is graded then $M^*$ is
graded, too. If $R$ is a graded $K$-algebra then $M$ is also a $K$-module
and the $K$-dual $M^{\vee}$ of $M$ is defined to be the graded module
$\Hom_K(M,K)$ where $K$ is considered as a graded module concentrated in
degree zero. Note that $R^{\vee}$ is the injective hull of $K^{\vee} \cong
K \cong R/{\mathfrak m}$ in the category of graded $R$-modules. If $\rankk
[M]_i < \infty$ for all integers $i$ then there is a canonical isomorphism
$M \cong M^{\vee \vee}$.

Now let $Z$ be a projective scheme over $K$. This means $Z = Proj (R)$
where $R$ is a graded $K$-algebra. For any sheaf $\cF$ on $Z$, we define
$H^i_*(Z,{\cF})=\bigoplus_{t\in \ZZ}H^i(Z,{\cF}(t))$.  In
this paper we will use ``vector bundle'' and ``locally free sheaf''
interchangeably.

Let $X$ be a non-empty projective subscheme of $Z$ with homogeneous
coordinate ring $A = R/I_X$. Then $I_X$ is a saturated ideal of $R$.
Recall that a homogeneous ideal $I$ in $R$ is {\em saturated} if $I =
\bigcup_{d \in {\ZZ}^+} [I:{\fm}^d ]$, where ${\fm} = (x_0 ,x_1
,\dots,x_n)$.   Equivalently, $I$ is saturated if and only if $I = H^0_*
(Z, {\cJ})$, where $\cJ$ is the sheafification of $I$.
\bigskip

\noindent {\it Generalized Koszul complexes}
\medskip

For more details with respect to the following discussion we refer to
\cite{BV} and \cite{Eisenbud-Buch}. The differences between these 
presentations and ours stem from the fact that we want to have all homomorphisms
graded (of degree  zero).

Let $R$ be a graded $K$-algebra and let $\ffi: F \to G$ be a homomorphism
of finitely generated graded $R$-modules. Then there are (generalized) Koszul
complexes $\cC_i(\ffi)$:
$$
0 \to \wedge^i F \otimes S_0(G) \to \wedge^{i-1} F \otimes S_1(G) \to
\ldots \to \wedge^0 F \otimes S_i(G) \to 0.
$$
Let $\cC_i(\ffi)^*$ be the $R$-dual of $\cC_i(\ffi)$. Suppose now that $F$ is a
free $R$-module of rank $f$. Then there are graded isomorphisms
$$
\wedge^f F \otimes (\wedge^j F)^* \cong \wedge^{f-j} F.
$$
Thus we can rewrite $\cC_i(\ffi)^* \otimes \wedge^f F$ as follows:
$$
0 \to \wedge^f F \otimes S_i(G)^* \to \wedge^{f-1} F \otimes
S_{i-1}(G)^* \to \ldots \to \wedge^{f-i} F \otimes S_0(G)^* \to 0.
$$
Note that $S_j(G)^*$ is the $j$th graded component of the divided power
algebra of $G^*$, but we won't need this fact.

Let's assume that also $G$ is a free $R$-module of rank, say, $g$ where $g
< f$. Then
$\ffi^*$ induces graded homomorphisms
$$
\nu_i: \wedge^{g+i} F \otimes \wedge^g G^* \to \wedge^i F.
$$
Put $r = f-g$. It turns out that for $i = 0,\ldots,r$ the complexes
$\cC_{r-i}(\ffi)^* \otimes \wedge^f F \otimes \wedge^g G^*$ and
$\cC_i(\ffi)$ can be spliced via $\nu_i$ to a complex $\cD_i(\ffi)$:
$$
0 \to \wedge^f F \otimes S_{r-i}(G)^* \otimes \wedge^g G^* \to \wedge^{f-1}
F \otimes
S_{r-i-1}(G)^* \otimes \wedge^g G^* \to \ldots
$$
$$
\to \wedge^{g+i} F
\otimes S_0(G)^* \otimes \wedge^g G^* \stackrel{\nu_i}{\longrightarrow}
\wedge^i F \otimes S_0(G) \to \wedge^{i-1} F \otimes S_1(G) \to
\ldots \to \wedge^0 F \otimes S_i(G) \to 0.
$$
The complex $\cD_0(\ffi)$ is called the Eagon-Northcott complex and
$\cD_1(\ffi)$ is called the Buchsbaum-Rim complex.

If we fix bases of $F$ and $G$ the map $\ffi$ can be described by a
matrix whose maximal minors generate an ideal which equals the image of
$\nu_0$. We denote this ideal by $I(\ffi)$. Its grade is at most
$f-g+1$. If $\ffi$ is general enough the complexes above have good
properties.

\begin{proposition} \label{EN-complexes_are_exact} Suppose $grade\, I(\ffi)
  = f-g+1$. Then it holds:
\begin{itemize}
\item[(a)] $\cD_i(\ffi)$ is acyclic where $i = 0,\ldots,f-g = r$.
\item[(b)] If $\ffi$ is a minimal homomorphism, i.e.\ $\im \ffi \subset \fm
  \cdot G$, then $\cD_0(\ffi)$ is a minimal free graded resolution of
  $R/I(\ffi)$
  and $\cD_i(\ffi)$ is a minimal free graded resolution of $S_i(\coker \ffi)$,
  $1 \leq i \leq r$.
\end{itemize}
\end{proposition}

The minimality of the resolutions in (b) follows by analyzing the maps
described above.
\bigskip

\noindent {\it Gorenstein rings and schemes}
\medskip

A graded $K$-algebra  $R$ is said to be Gorenstein if it has finite
injective dimension (cf.\ \cite{Bruns-Herzog}, Definition 3.1.18). Over a
Gorenstein ring duality theory is particularly
simple. We denote the index of regularity of a graded ring by $r(R)$. If
$R$ is just the polynomial ring $K[x_0, \dots , x_n]$ then $r(R) = -n$.  We
will use the following duality result (cf., for example, \cite{SV2},
Theorem 0.4.14).

\begin{lemma} \label{duality}  Let $M$ be a graded $R$-module where $R$ is
  a Gorenstein ring of dimension $n$. Then we have  for all $i \in
  \mathbb{Z}$ natural isomorphisms of graded $R$-modules
$$
H^i_{{\mathfrak m}} (M)^{\vee} \cong \Ext^{n-i}_R(M,R)(r(R)-1).
$$
\end{lemma}

Let $M$ be a graded  $R$-module where $n = \dim R$ and $d = \dim M$. Then
$$
K_M = \Ext^{n-d}_R(M,R) (r(R)-1)
$$
is said to be the {\itshape canonical module} of $M$. Usually the canonical
module is defined as the module representing the functor $\HH^d(M \otimes_R
\__{})^{\vee}$ if such a module exists. If $R$ is Gorenstein it does and is
just the module defined above (cf.\ \cite{S}).

We say that $M$ has cohomology of finite length if the cohomology modules
$\HH^i(M)$ have finite length for all $i < \dim M$. It is well-known that
$M$ has cohomology of finite length if and only if $M$ is equidimensional
and locally Cohen-Macaulay.

Let now $Z = Proj (R)$ be a  projective scheme over $K$. Then $Z$ is
said to be {\itshape arithmetically Gorenstein} and {\itshape
  arithmetically Cohen-Macaulay}
respectively if the homogeneous coordinate ring $R$ of $Z$ is Gorenstein and
Cohen-Macaulay respectively.

For a closed subscheme $X$ of $Z$, with homogeneous coordinate ring $A =
R/I_X$ we will refer to the canonical module of $A$ also as the canonical
module of $X$. Moreover, we say that $X$ has finite projective dimension if
$A$
has finite projective dimension as an $R$-module.

Assume that $Z$ is arithmetically Gorenstein. One of the things we  shall
be interested in is to describe when certain subschemes $X$ of $Z$
are arithmetically Gorenstein, too.   To do this, it
is enough to show that $X$ is arithmetically Cohen-Macaulay, with
Cohen-Macaulay type  $1$ provided
$X$ has finite projective dimension. In this case $X$ is
defined by a Gorenstein ideal $I = I_X \subset R$.

Recall that the {\em Cohen-Macaulay type} of an arithmetically Cohen-Macaulay
projective scheme $X$ with finite projective dimension can be defined to be
the rank of the last free module
occurring in a minimal free resolution of the saturated ideal of $X$. It
is equal to  the number of minimal generators of the canonical module of
$X$.


\section{Buchsbaum-Rim sheaves}

 From now on we will always assume that $Z$ is a projective arithmetically
Gorenstein scheme over the field $K$.  We denote its  dimension by $n$ and
its  homogeneous coordinate ring by $R$.

Let $\cF$ and $\cG$ be locally free sheaves of ranks $f$ and $g$
respectively on $Z$.  Let $\ffi: \cF \to \cG$ be a generically surjective
morphism. Since the
construction of the generalized Koszul complexes as described in the
previous section globalizes,  we can
associate to $\ffi$ several complexes. The most familiar are the
Eagon-Northcott complex
$$
0 \to \wedge^f \cF \otimes S_{f-g}(\cG)^* \otimes \wedge^g \cG^* \to
\wedge^{f-1} \cF \otimes
S_{f-g-1}(\cG)^* \otimes \wedge^g \cG^* \to \ldots
$$
$$
\to \wedge^{g} \cF
\otimes S_0(\cG)^* \otimes \wedge^g \cG^* \stackrel{\wedge^g
  \ffi}{\longrightarrow} \cO_Z \to 0
$$
and the Buchsbaum-Rim complex
$$
0 \to \wedge^f \cF \otimes S_{f-g-1}(\cG)^* \otimes \wedge^g \cG^* \to
\wedge^{f-1} \cF \otimes
S_{f-g-2}(\cG)^* \otimes \wedge^g \cG^* \to \ldots
$$
$$
\to \wedge^{g+1} \cF
\otimes S_0(\cG)^* \otimes \wedge^g \cG^* \to \cF
\stackrel{\ffi}{\longrightarrow} \cG  \to 0.
$$
Moreover, Proposition~\ref{EN-complexes_are_exact} implies that these
complexes are acyclic if the degeneracy locus of $\ffi$ has the expected
codimension $f-g+1$ in $Z$. This lead us to the following definition.

\begin{definition} With the notation above suppose that the degeneracy locus of
  $\ffi$ has codimension $f-g+1$. Then we call the cokernel of the map
  between $\wedge^{g+1+i} \cF \otimes S_i(\cG)^* \otimes \wedge^g
  \cG^*$ and  $\wedge^{g+i} \cF \otimes S_{i-1}(\cG)^* \otimes
  \wedge^g \cG^*$ an {\it $i^{th}$ local Buchsbaum-Rim sheaf} ($1 \leq i
  \leq f-g$) and denote it by $\cB^{\ffi}_i$.
\end{definition}

Note that the $i^{th}$  local Buchsbaum-Rim sheaf associated to $\ffi$ is just
the $(i+1)^{st}$ syzygy sheaf of $\coker \ffi$.

The following result is a generalization of Proposition 2.10 of
\cite{KMNP}. Thanks to our set-up the proof given there works here, too.

\begin{proposition} Let $\cB$ be a first local Buchsbaum-Rim sheaf
  associated to a morphism $\ffi$. Let $X$ denote  the degeneracy locus of
  $\ffi$. Let $S$ be the zero-locus of a section $s \in H^0(Z,\cB)$ and let
  $T$ be the zero-locus of a section $t \in H^0(Z,\cB^*)$. Then it holds $X
  \subset S$ and $X \subset T$.
\end{proposition}

Now we put stronger assumptions on the sheaves $\cF$ and $\cG$.

\begin{definition} \label{BR-sheaf} Suppose in addition that the modules $F =
  H^0_*(Z,\cF)$ and $G = H^0_*(Z,\cG)$ are  free
  $R$-modules. Then the sheaf $\cB^{\ffi}_i$ is called an {\it $i^{th}$
    Buchsbaum-Rim sheaf}. For simplicity a first Buchsbaum-Rim sheaf is
  just called a {\it Buchsbaum-Rim sheaf} and denoted by $\cBf$.

By abuse of notation we denote the homomorphism $F \to G$ induced by $\ffi$
again by $\ffi$. Moreover, we put $\Bf = H^0_*(Z,\cBf)$ and $\Mf = \coker
\ffi$, so that we have an exact sequence
$$
0 \to B_{\ffi} \to F \stackrel{\ffi}{\longrightarrow} G \to  M_{\ffi} \to 0.
$$
We call $\Bf$ a {\it Buchsbaum-Rim module}.
\end{definition}

\begin{remark} In this paper we will consider some degeneracy loci
  associated to  Buchsbaum-Rim sheaves. These investigations were motivated
  by the work of the first and the third author in \cite{Mig-P_gorenstein}.

  Note that in \cite{KMNP}
  zero-loci of regular sections of the {\it dual} of a Buchsbaum-Rim sheaf
  over projective space have been characterized as determinantal subschemes
  which are generically complete intersections.

If $g=1$ then the sheaf $\cB^{\ffi}_{f-1}$ is the dual of a Buchsbaum-Rim
sheaf. Thus it seems to be rewarding to study higher Buchsbaum-Rim sheaves,
too.
\end{remark}

\begin{remark} \label{properties_of_BR-sheaf} (i) With the notation above
  the Buchsbaum-Rim sheaf $\cBf$ has rank $r = f - g$. Our assumptions
  imply that $r \leq n = \dim Z$. Moreover, $\cBf$ is locally free if and
only if $n
  = r$.

(ii) As a second syzygy sheaf over an arithmetically Gorenstein scheme a
Buchsbaum-Rim sheaf $\cBf$ is reflexive, i.e., the natural map $\cBf
\otimes \cBf^* \to \cO_Z$ induces an isomorphism $\cBf \cong \cBf^{**}$.

Similarly, a Buchsbaum-Rim module is a reflexive $R$-module.
\end{remark}

The following result will become important later on.

\begin{lemma} \label{prop-of-exteriour-powers} With the above notation let
  $\Bf$ be a Buchsbaum-Rim module of rank $r$. Then it holds:
\begin{itemize}
\item[(a)] For $i = 1,\ldots r$ the
  module $\wedge^i \Bf^*$ is a $(r-i+1)$-syzygy of the perfect
  module $S_{r-i}(\Mf) \otimes \wedge^f F^* \otimes \wedge^g G$ and
  is resolved by $\cC_i(\ffi)^*$.
\item[(b)] For $i = 0,\ldots r$ there are isomorphisms
$$
\wedge^{r-i} \cBf^* \otimes \wedge^f F \otimes \wedge^g G^* \cong (\wedge^i
\cBf^*)^*.
$$
\item[(c)] For $i < r$ we have
$$
\HH^j(\wedge^i \Bf^*) \cong \left \{ \begin{array}{ll}
0 & \mif j \leq n \; \mbox{and} \; j \neq n+1-i \\
S_{i}(\Mf)^{\vee}(1 - r(R)) & \mif j = n+1-i.
\end{array} \right.
$$
\end{itemize}
\end{lemma}

\begin{proof} The assumption $\codim I(\ffi) = r+1$ ensures that
  $\wedge^i \Bf^*$ is resolved by $\cC_i(\ffi)^*$ (cf.\ \cite{BV},
  Remark 2.19).  Hence the first claim follows by
Proposition~\ref{EN-complexes_are_exact}.

The latter result also implies the second claim. The isomorphism is induced
by the map $\nu_i$ (cf.\ \cite{BV}, Remark 2.19).

In order to prove the third claim we observe that by Lemma~\ref{duality}
$$
\HH^{n+1-j}(\wedge^i \Bf^*) \cong \Ext^j(\wedge^i
\Bf^*,R)^{\vee}(1-r(R)).
$$
But we know already that $\Ext^j(\wedge^i \Bf^*,R)$ can be computed as the
$(i-j)^{th}$ homology module of the complex $\cC_{i}(\ffi)$ which is part
of the acyclic complex $\cD_{i}(\ffi)$. Now our assertion follows.
\end{proof}

\begin{remark} \label{Bott-formula} (i) The previous result implies for the
  Buchsbaum-Rim sheaf $\cBf$ that it is just the sheafification
  $\widetilde{\Bf}$ of the Buchsbaum-Rim module $\Bf$ and
$$
H_*^j(\wedge^i \cBf^*) \cong \left \{ \begin{array}{ll}
0 & \mif 1 \leq j <  n \; \mbox{and} \; j \neq n-i \\
S_{i}(\Mf)^{\vee}(1 - r(R)) & \mif j = n-i.
\end{array} \right.
$$
(ii) Let us consider the example where $R = K[x_0,\ldots,x_n]$ is the
polynomial ring, $F = R(-1)^{n+1}, G = R$ and $\ffi: F \to G$ a general
map. Then $\Mf \cong K$ and $\Bf = \ker \ffi$ is a Buchsbaum-Rim module and
the corresponding
Buchsbaum-Rim sheaf is just the cotangent bundle on $\PP^n$. Via Serre
duality we see that in this case Lemma~\ref{prop-of-exteriour-powers}(c) is
just the dual version of the Bott formula for the cohomology of $\cOP^j =
\wedge^j \cBf$.
\end{remark}

According to our explicit description of the complexes $\cD_{i}(\ffi)$,
Proposition~\ref{EN-complexes_are_exact} implies.

\begin{lemma} \label{can-module-of-symm-powers}  For $i = 1,\ldots r-1$ there
  are isomorphisms
$$
\Ext^{r+1}(S_i(M_{\ffi}),R) \cong S_{r-i}(M_{\ffi}) \otimes \wedge^f F^*
\otimes \wedge^g G.
$$
\end{lemma}

Finally, we want to derive a cohomological characterization of
Buchsbaum-Rim sheaves. It turns out that they are particular
Eilenberg-MacLane sheaves. Recall that an $R$-module $E$ is called an {\it
  Eilenberg-MacLane module} of depth $t$, $0 \leq t \leq n+1$ if
$$
\HH^j(E) = 0 \quad \mbox{for all} \;  j \neq t \; \mbox{where} \; 0 \leq j
\leq n.
$$
Similarly, a sheaf $\cE$ on $Z$ is said to be an {\it Eilenberg-MacLane
  sheaf} if $H^0_*(\cE)$ is an Eilenberg-MacLane module.

We will need the following result which is shown in \cite{habil} as
Theorem I.3.9.

\begin{lemma} \label{Eilenberg-MacLane}  Let $E$ be a
  reflexive module of depth $t \leq n$. Then $E$ is an Eilenberg-MacLane
  module with finite projective dimension if and only if $E^*$ is an
  $(n+2-t)$-syzygy of a module $M$ of dimension $\leq t-2$. In this case it
  holds
$$
M \cong \HH^t(E)^{\vee}(1-r(R)).
$$
\end{lemma}

Now we are ready for our cohomological description of Buchsbaum-Rim
sheaves.

\begin{proposition} \label{BR-sheaf-characterization} A sheaf $\cE$ on $Z$
  is a Buchsbaum-Rim sheaf if and only if $E = H^0_*(\cE)$ is a reflexive
  Eilenberg-MacLane module with finite projective dimension and rank $r
  \leq n$ such that $\HH^{n-r+2}(E)^{\vee}$ is a perfect $R$-module of
  dimension  $n-r$ if $r \geq 2$.
\end{proposition}

\begin{proof} First let us assume that $\cE$ is a Buchsbaum-Rim sheaf. Then
  we have by definition for $E = H^0_*(\cE)$ that it has a rank, say $r$,
  and sits in an exact sequence
$$
0 \to E \to F \stackrel{\ffi}{\longrightarrow} G \to  M_{\ffi} \to 0
$$
where $F$ and $G$ are free modules and $I(\ffi)$ has the expected
codimension $r+1$. Due to Remark~\ref{properties_of_BR-sheaf} $E$ is
reflexive. Furthermore, $E$ has finite projective dimension since $\Mf$
does by Lemma~\ref{prop-of-exteriour-powers}. If $r=1$ it follows that $E$
is just $R(m)$ for some integer $m$.

Let $r \geq 2$. Then the exact sequence above and the Cohen-Macaulayness of
$\Mf$ imply that $E$ is an Eilenberg-Maclane module of depth $n-r+2$ and
$$
\HH^{n-r+2}(E)^{\vee} \cong \HH^{n-r}(\Mf)^{\vee} \cong
\Ext^{r+1}(\Mf,R)(r(R)-1) \cong S_{r-1}(\Mf) \otimes \wedge^f F^* \otimes
\wedge^g G(r(R) - 1)
$$
where the latter isomorphisms are due to Lemma~\ref{duality} and
Lemma~\ref{can-module-of-symm-powers}. Since $S_{r-1}(\Mf)$ is a perfect module
of dimension $n-r$ by Lemma~\ref{prop-of-exteriour-powers} we have shown
that the
conditions in the statement are necessary.

Now we want to show sufficiency. Since a reflexive module of rank $1$ with
finite projective dimension must be free we are done if $r=1$. Now let $r
\geq 2$.  By assumption $M = \HH^{n-r+2}(E)^{\vee}$
is a perfect module of dimension $n-r$. Since $E^*$ is an $r$-syzygy of $M$
due to Lemma~\ref{Eilenberg-MacLane} we obtain that $E^*$ is a reflexive
Eilenberg-MacLane module of depth $n$ with finite projective dimension and
$$
\HH^n(E^*) \cong \HH^{n-r}(M).
$$
Therefore Lemma~\ref{Eilenberg-MacLane} applies also to $E^*$ and says that
$E^{**} \cong E$ is a $2$-syzygy of $\HH^{n-r}(M)^{\vee}(1-r(R))$. This
means that there is an exact sequence
$$
0 \to E \to F \stackrel{\ffi}{\longrightarrow} G \to K_M(1-r(R)) \to 0
$$
where $F$ and $G$ are free modules with $\rank F = \rank G + r$. Since
$\Rad I(\ffi) = \Rad \Ann_R K_M$ and $K_M$ has dimension $n-r$ it follows
that $I(\ffi)$ has the maximal codimension $r+1$. Thus $E$ is a Buchsbaum-Rim
module completing the proof.
\end{proof}

Since any module over a regular ring has finite projective dimension the
last result takes a simpler form for sheaves on $\PP^n$.

\begin{corollary} \label{Bu-Rim-sheafs-on-proj-space} A sheaf $\cE$ on
  $\PP^n$ is a Buchsbaum-Rim sheaf if and only if $\cE$ is an reflexive
  Eilenberg-MacLane sheaf of rank $r \leq n$ such that $H^i_*(\PP^n,\cE) =
  0$ if $i \neq 0, n-r+1, n+1$ and $H^{n-r+1}_*(\cE)^{\vee}$ is a
  Cohen-Macaulay module of dimension $n-r$.
\end{corollary}

  From this result we see again that the cotangent bundle on projective space
is a Buchsbaum-Rim sheaf.


\section{The cohomology of the degeneracy loci}

Consider a morphism $\psi: \cP \to \cBf$ of sheaves of rank $t$ and
$r$ respectively on the arithmetically Gorenstein scheme $Z = Proj (R)$
where $\cBf$ is a Buchsbaum-Rim sheaf and $H^0_*(Z, \cP)$ is a free
$R$-module. If $t = 1$ then $\psi$ is just a section of some twist of
$\cBf$. Thus we refer to $\psi$ as multiple sections of $\cBf$.

Throughout this paper we suppose that the ground field $K$ is infinite and
that the degeneracy locus $S$ of $\psi$ 
has (the expected) codimension $r - t + 1 \geq 2$ in $Z$. If $t=1$ then $S$
is just the zero-locus of a regular section of a Buchsbaum-Rim sheaf.

It will turn out that $S$ is often not equidimensional. Thus we are also
interested in the top-dimensional part $X$ of $S$, i.e. the union of the
highest-dimensional components of $S$. The aim of this section is to
compute the cohomology modules of $S$ and $X$ respectively. An
Eagon-Northcott complex involving $\cBf$ will play an essential
role. Observe that in contrast to the situation in the previous section
where $\cF$ and $\cG$ were locally free the
sheaf $\cBf$ is in general not locally free.

Our approach will be algebraic. Taking global sections we obtain an
$R$-homomorphism
$$
\psi : P \to \Bf
$$
where $P$ is a free $R$-module of rank $t, 1 \leq t < r$. The first aim is
to derive the complex mentioned above. We follow the approach described in
Section~\ref{preliminaries}. The $R$-dual of the Koszul complex
$\mathcal{C}_{r-t}(\psi^*)$ is
$$
0 \to (\wedge^0 \Bf^* \otimes S_{r-t}(P^*))^* \to (\Bf^* \otimes
S_{r-t-1}(P^*))^* \to \ldots \to
(\wedge^{r-t-1} \Bf^* \otimes P^*)^* \to (\wedge^{r-t} \Bf^* \otimes
S_0(P^*))^*.
$$
Using the isomorphisms in Lemma~\ref{prop-of-exteriour-powers} and
$S_j(P^*)^* \cong S_j(P)$ we can
rewrite $\mathcal{C}_{r-t}(\psi^*)^* \otimes \wedge^f F^* \otimes \wedge^g G$
as follows:
$$
0 \to \wedge^r \Bf^* \otimes S_{r-t}(P) \to \wedge^{r-1} \Bf^* \otimes
S_{r-t-1}(P) \to \ldots \to
\wedge^{t+1} \Bf^* \otimes P \to \wedge^{t} \Bf^* \otimes S_0(P).
$$
The image of the map $\wedge^{t} \psi^* : \wedge^t \Bf^* \to
\wedge^t P^*$ is (up to degree shift) an ideal of $R$ which we
denote by $I(\psi)$ or just $I$, i.e.\ $\im \wedge^{t} \psi^* = I \otimes
\wedge^t P^*$. Note that the saturation of $I$ is the
homogeneous ideal $I_S$ defining the degeneracy locus $S$. Thus, using
$\wedge^{t} \psi^*$ we can continue the complex above on the right-hand
side and obtain the desired Eagon-Northcott complex
\begin{eqnarray}
\lefteqn{E_{\bullet} : \quad 0 \to \wedge^r \Bf^* \otimes S_{r-t}(P) \to
  \wedge^{r-1} \Bf^* \otimes
S_{r-t-1}(P) \to \ldots} \\
& &  \to \wedge^{t+1} \Bf^* \otimes P \to \wedge^{t} \Bf^* \to I \otimes
\wedge^t P^* \to 0.  \nonumber
\end{eqnarray}

The next result shows that the first cohomology modules of $R/I(\psi)$
vanish.

\begin{lemma} \label{depth-estimate} The depth of $R/I$ is at least
  $n-r$.
\end{lemma}

\begin{proof} Let $r=2$. Then $\Bf$ has depth $n$ and $t=1$. Thus we have
  an exact sequence
$$
0 \to R(a) \stackrel{\psi}{\longrightarrow} \Bf \to I(b) \to 0
$$
where $a, b$ are integers. It provides the claim.

Now let $r \geq 3$.  We choose a
  sufficiently general linear form $l \in R$. For short we denote the
  functor $\_\_\otimes_R \overline{R}$ by  $^-$ where $\overline{R} = R/l R$.
Let $\alpha$ be the map
  $\Hom_{\overline{R}}(\overline{\psi},\overline{R}) :
    \Hom_{\overline{R}}(\overline{\Bf},\overline{R}) \to
    \Hom_{\overline{R}}(\overline{P},\overline{R})$ and define the
    homogeneous ideal $J \subset \overline{R}$ by  $J \otimes \wedge^t
    \overline{P}^* = \im \wedge_{\overline{R}}^t \alpha$. Our first claim
    is that $\overline{I} = J$ provided $n > r$.

Multiplication by $l$ provides the commutative diagram
$$
\begin{array}{ccccccc}
0 \to & P(-1) & \to & P  & \to & \overline{P} & \to 0 \\
& \downarrow\!{\scriptstyle \psi} & & \downarrow\!{\scriptstyle \psi} & &
\downarrow\!{\scriptstyle \overline{\psi}}  &\\
0 \to & \Bf(-1) & \to & \Bf  & \to & \overline{\Bf} & \to 0 \\
\end{array}
$$
Dualizing gives the commutative diagram
$$
\begin{array}{ccccccc}
0 \to & \Bf^* & \to & \Bf^*(1)  & \to & \ER^1(\overline{\Bf},R) & \to
\ER^1(\Bf,R) = 0 \\
& \downarrow\!{\scriptstyle \psi^*} & & \downarrow\!{\scriptstyle \psi^*} & &
\downarrow\!{\scriptstyle \beta}  &\\
0 \to & P^* & \to & P^*(1)  & \to & \ER^1(\overline{P},R) & \to \ER^1(P,R)
= 0 \\
\end{array}
$$
where the vanishings on the right-hand side are due to duality and the fact
that $\Bf$ is an Eilenberg-MacLane module of depth  $n-r+2 \neq n$.
Using $\ER^1(\overline{\Bf},R)(-1) \cong
\Hom_{\overline{R}}(\overline{\Bf},\overline{R})$ and
$\ER^1(\overline{P},R)(-1) \cong
\Hom_{\overline{R}}(\overline{P},\overline{R})$ (cf., for example,
\cite{Bruns-Herzog}, Lemma 3.1.16) we see
that $\beta$ can be identified with $\alpha$ as well as $\psi^* \otimes
\overline{R}$. It follows
$$
J \otimes \wedge^t {\overline{P}}^* = \im \wedge_{\overline{R}}^t \alpha = \im
\wedge_{\overline{R}}^t
(\psi^*
  \otimes \overline{R}) \cong \im (\wedge_R^t \psi^*
  \otimes \overline{R}) \cong (\im \wedge_R^t \psi^*)
  \otimes \overline{R} = I \otimes \wedge^t P^* \otimes \overline{R}
$$
and thus $J = \overline{I}$ as we wanted to show.

The second commutative diagram also provides $\Bf^* \otimes \overline{R}
\cong \Hom_{\overline{R}}(\overline{\Bf},\overline{R})$. It follows
$$
\wedge^t_{\overline{R}} \Hom_{\overline{R}}(\overline{\Bf},\overline{R})
\cong \wedge^t_{\overline{R}} (\Bf^* \otimes \overline{R}) \cong
(\wedge^t_R \Bf^*) \otimes \overline{R}.
$$
Thus we have an exact commutative diagram
$$
\begin{array}{ccccccc}
& 0 & & 0 & & 0 &\\
& \downarrow & & \downarrow & & \downarrow  &\\
0 \to & C(-1) & \to & \wedge^t \Bf^*(-1)  & \to & I(-1) \otimes \wedge^t
P^* & \to 0 \\
& \downarrow\!{\scriptstyle l} & & \downarrow\!{\scriptstyle l} & &
\downarrow\!{\scriptstyle l}  &\\
0 \to & C & \to & \wedge^t \Bf^*  & \to & I \otimes \wedge^t P^* & \to 0 \\
&  & & \downarrow & &  &\\
0 \to & L & \to & \wedge^t_{\overline{R}}
\Hom_{\overline{R}}(\overline{\Bf},\overline{R})   & \to & J \otimes \wedge^t
{\overline{P}}^* & \to 0 \\
&  & & \downarrow & &  &\\
&  & & 0 & &  &
\end{array}
$$
where $C = \ker \wedge_R^t \psi^*$ and $L = \ker
\wedge_{\overline{R}}^t \alpha$. Since $\wedge^t \Bf^*$ has depth $n+1-t
   \geq n+2-r$ due to Lemma~\ref{prop-of-exteriour-powers} the first row
   shows that our assertion is equivalent to
   $\depth C \geq n+2-r$. In order to show this we induct on $n-r$. If $n =
   r$ the claim follows by the exact cohomology sequence induced by  the
   top-line of the
   previous diagram and $\depth \wedge^t \Bf^* \geq 2$. Let $n > r$. Then
   our first claim and the Snake
   lemma applied to the diagram above imply $\overline{C} \cong L$. The
   induction hypothesis applies to $L$ and we obtain
$$
0 < n+1-r \leq \depth L < \depth C
$$
completing the proof.
\end{proof}

For computing the other cohomology modules of $R/I(\psi)$ we use the
Eagon-Northcott complex above. In order to ease notation let us  write
$E_{\bullet}$ as
$$
E_{\bullet}  \colon \quad 0 \to E_r \stackrel{\delta_r}{\longrightarrow} E_{r-1}
\stackrel{\delta_{r-1}}{\longrightarrow} \cdots \to E_t
\stackrel{\delta_t}{\longrightarrow} I(\psi) (p)  \to 0
$$
where
$$
E_i = \wedge^i \Bf^* \otimes S_{i - t}(P) \quad \mbox{and} \quad R(p)
\cong \wedge^t P^*.
$$
The number of minimal generators of an $R$-module $N$ is denoted by
$\mu(N)$.

\begin{proposition} \label{cohomology-of locus} Put $I = I(\psi)$. Then we
have:
\begin{itemize}
\item[(a)] For  $j\neq \dim R/I = n+t-r$ it holds
$$
H^j_m (R/I) \cong \left\{ \begin{array}{ll}
S_i(M_{\ffi})^\vee \otimes
S_{i-t} (P) \otimes \wedge^t P (1-r(R)) & \mbox{if } j=n+t-2i \quad
\mbox{where }
\max \{t, \frac{r+1}{2}\} \leq i \leq \left \lfloor \frac{r+t}{2} \right
\rfloor  \\
0 & \mbox{otherwise}.
\end{array} \right.
$$
\item[(b)] The canonical module satisfies
$$
\mu(K_{R/I}) \leq \binom{r-1}{t-1} + \left\{
\begin{array}{ll}
\binom{\frac{r}{2}+g-1}{g-1} \cdot \binom{\frac{r}{2} - 1}{t-1} & \mbox{if
  $r$ is even and } 1\leq t \leq \frac{r}{2} \\
0 & \mbox{otherwise}.
\end{array}\right.
$$
\end{itemize}
\end{proposition}

\begin{proof} We consider the Eagon-Northcott complex $E_{\bullet}$ above.
According to the Lemma~\ref{prop-of-exteriour-powers} $\Bf^*$ is an $r$-syzygy.
Therefore $\Bf^*$ is locally free in codimension $r$. It follows that  the
Eagon-Northcott complex $E_{\bullet}$ is exact in codimension $r$.
Therefore its
homology modules $H_i(E_{\bullet})$ have dimension $\leq n-r$. Thus the exact
sequence
$$
0 \to \im \delta_{i+1} \to \ker \delta_i \to H_i(E_{\bullet}) \to 0
$$
implies
$$
\HH^j(\ker \delta_i) \cong \HH^j(\im \delta_{i+1}) \quad \mif j \geq
n-r+2. \leqno(1)
$$
Moreover, there are exact sequences
$$
0 \to \ker \delta_i \to E_i \to \im \delta_i \to 0
$$
inducing  exact sequences
$$
\HH^j(E_i) \to \HH^j(\im \delta_i) \to \HH^{j+1}(\ker \delta_i)
\to \HH^{j+1}(E_i) \leqno(2)
$$
where the injectivity or surjectivity of the map in the middle can be
checked by means of Lemma~\ref{prop-of-exteriour-powers}. This map is an
isomorphism if $j \neq n-i, n+1-i, n, n+1$.

Let us consider the map $\delta_r$. Due to our assumption the map $\psi^*:
\Bf^* \to P^*$ is generically surjective. Thus the same applies to the
Koszul map $\Bf^* \otimes S_{r-t-1}(P^*) \to S_{r-t}(P^*)$ which is induced
by $\psi$. It follows that
the $R$-dual of this map is injective. But the latter is (up to a degree
shift) just $\delta_r$. Hence we have seen that $\im \delta_r \cong E_r$.

According to Lemma~\ref{depth-estimate} it suffices to consider $\HH^j(R/I)$
where $n-r \leq j \leq n+t-r = \dim R/I$. For this we distinguish several
cases. \\
{\it Case 1}: Let us assume that $n-t < j \leq n+t-r$. This can occur if
and only if $t \geq \frac{r+1}{2}$.

Using (2) and (1) in alternating order we get
\begin{eqnarray*}
\lefteqn{\HH^j(R/I)(p) \cong \HH^{j+1}(\im \delta_t) \cong \HH^{j+2}(\ker
  \delta_t) \cong \HH^{j+2}(\im \delta_{t+1}) \cong \ldots} \\
& & \cong \HH^{r+j-t}(\im \delta_{r-1}) \hookrightarrow \HH^{r+j+1-t}(\ker
\delta_{r-1}) \cong \HH^{r+j+1-t}(\im \delta_{r}) \cong \HH^{r+j+1-t}(E_r)
\end{eqnarray*}
where the injection holds true because we have by our assumptions $n+3-r
\leq n+2-t \leq j+1 \leq r+j-t \leq n$.

It follows by Lemma~\ref{prop-of-exteriour-powers} and Lemma~\ref{duality}
that
$$
\HH^j(R/I) = 0 \quad \mif n-t < j < n+t-r
$$
and in case $j = \dim R/I$ for the canonical module
$$
\mu(K_{R/I}) \leq \mu(\HH^{n+1}(E_r)^\vee) = \mu(E_r) = \mu(S_{r-t}(P)) =
\binom{r-1}{t-1}.
$$
{\it Case 2}: Let us assume that $n-r \leq j= n+t-1-2i \leq \min \{n+t-r,
n-t \}$, i.e.\ $\max \{t, \frac{r-1}{2} \} \leq i \leq \frac{r+t-1}{2}$.

Using (2) and (1) again we obtain
\begin{eqnarray*}
\lefteqn{\HH^{n+t-1-2i}(R/I)(p) \cong \HH^{n+t-2i}(\im \delta_t) \cong
  \HH^{n+t+1-2i}(\ker \delta_t) \cong \HH^{n+t+1-2i}(\im \delta_{t+1})
  \cong \ldots} \\
& \cong & \HH^{n-i}(\im \delta_i) \hookrightarrow \HH^{n+1-i}(\ker
\delta_i) \cong \ldots \\
& \cong & \HH^{n+r-1-2i}(\im \delta_{r-1}) \hookrightarrow
\HH^{n+r-2i}(\ker \delta_{r-1}) \cong \HH^{n+r-2i}(\im \delta_{r}) \cong
\HH^{n+r-2i}(E_r).
\end{eqnarray*}
Since $E_r$ is free we get
$$
\HH^{n+t-1-2i}(R/I) = 0 \quad \mif i \geq \frac{r}{2}
$$
and in case $i = \frac{r-1}{2}$ we obtain the same bound for the number of
minimal generators of the canonical module as in Case 1. \\
{\it Case 3}: Let us assume that $n-r \leq j= n+t-2i \leq \min \{n+t-r,
n-t \}$, i.e.\ $\max \{t, \frac{r}{2} \} \leq i \leq \frac{r+t}{2}$.

If follows similarly as before
\begin{eqnarray*}
\lefteqn{\HH^{n+t-2i}(R/I)(p) \cong \HH^{n+t+1-2i}(\im \delta_t)} \\
&  \cong &
  \HH^{n+t+2-2i}(\ker \delta_t) \cong \HH^{n+t+2-2i}(\im \delta_{t+1})
  \cong \ldots \\
& \cong & \HH^{n+1-i}(\im \delta_i).
\end{eqnarray*}
Now we look at the exact sequence
$$
\HH^{n+1-i}(\ker \delta_i) \to \HH^{n+1-i}(E_i) \to \HH^{n+1-i}(\im
\delta_i) \to \HH^{n+2-i}(\ker \delta_i). \leqno(3)
$$
We use again (1) and (2) in order to obtain information on the modules on
the left-hand and on the right-hand side. This provides
$$
\HH^{n+1-i}(\ker \delta_i) \cong \HH^{n+1-i}(\im \delta_{i+1}) \cong \ldots
\cong \HH^{n+r-2i}(\im \delta_r) \cong \HH^{n+r-2i}(E_r) = 0
$$
because $n+r-2i \leq n$ and
\begin{eqnarray*}
\lefteqn{\HH^{n+2-i}(\ker \delta_i) \cong \HH^{n+2-i}(\im \delta_{i+1})
  \cong \ldots } \\
& & \cong \HH^{n+r-2i}(\im \delta_{r-1}) \hookrightarrow \HH^{n+r-2i}(\ker
\delta_r) \cong \HH^{n+r+1-2i}(\im \delta_r) \cong \HH^{n+r+1-2i}(E_r) 
\end{eqnarray*}
where the last module vanishes if and only if $i \neq \frac{r}{2}$.

Therefore (3) yields if $i \neq \frac{r}{2}$
$$
\HH^{n+t-2i}(R/I)(p) \cong \HH^{n+1-i}(E_i) \cong
\HH^{n+1-i}(\wedge^i\Bf^*) \otimes S_{{\frac{r}{2}}-t}(P) \cong
S_{\frac{r}{2}}(M_{\ffi})^{\vee} 
\otimes S_{{\frac{r}{2}}-t}(P) (1 - r(R)).
$$
In case $i = \frac{r}{2}$ taking $K$-duals of (3) furnishes the exact
sequence
$$
E_r^* \to \Ext^{r-t+1}(R/I,R)(-p) \to S_{\frac{r}{2}}(M_{\ffi}) \otimes
S_{{\frac{r}{2}}-t}(P). 
$$
It follows for the canonical module
$$
\begin{array}{lcl}
\mu(K_{R/I}) & \leq & \mu(S_{\frac{r}{2}} (M_{\ffi})) \otimes S_{\frac{r}{2}-t}
(P) + \mu(E_r) \\
& = & \binom{\frac{r}{2}+g-1}{g-1} \cdot \binom{\frac{r}{2}-1}{t-1} +
\binom{r-1}{t-1}.
\end{array}
$$
Our assertions are now a consequence of the results in the three cases
above.
\end{proof}

\begin{corollary} \label{depth-formula} For the depth of the coordinate
  ring we have
$$
\depth R/I(\psi) = \left \{ \begin{array}{ll}
n-r & \mif r+t \; \mbox{is even} \\
n-r+1 & \mif r+t \; \mbox{is odd}.
\end{array} \right.
$$
\end{corollary}

Put $e = \depth R/I(\psi)$. Then the only non-vanishing cohomology modules
of $R/I(\psi)$ besides $\HH^{n+t-r}(R/I(\psi))$ are $\HH^{e+2k}(R/I(\psi))$
where $k$ is an integer with $0 \leq k \leq \frac{1}{2} [\min \{
n-t,n+t-r-1 \} - e]$.
\smallskip

It has already been observed in \cite{Mig-P_gorenstein} that $I(\psi)$ is
not always an unmixed ideal. This gives rise to consider the ideal
$J(\psi)$ which is defined as the intersection of the primary components of
$I(\psi)$ having maximal dimension. We denote by $X$ the subscheme of $Z$
defined by $J(\psi)$ and call it the top-dimensional part of the degeneracy
locus $S$.

Our next aim is to clarify the relationship between $S$ and $X$. For this
we need a cohomological criterion for unmixedness stated as  Lemma
III.2.3 in \cite{habil}.

\begin{lemma} \label{unmixedness-criterion} Let $I \subset R$ be a
  homogeneous ideal. Then $I$ is unmixed if and only if
$$
\dim \Supp(\HH^i(R/I)) < i \quad \mbox{for all} \; i < \dim R/I
$$
where we put $\dim \Supp(M) = - \infty$ if $M = 0$.
\end{lemma}

Now we can show.

\begin{proposition} \label{coho_of_top-dimensional} Let $I = I(\psi)$ and
  $J = J(\psi)$. Then it holds:
\begin{itemize}
\item[(a)]  $I$ is unmixed if and only if $r+t$ is odd.
\item[(b)] If $r+t$ is even then $I$ has a primary component of codimension
  $r+1$. Let $Q$ be the intersection of all those components. Then we have $I =
  J \cap Q$ and
$$
H^j_m(R/J) \cong \left\{ \begin{array}{ll}
H^j_m (R/I) & \mbox{if } j\neq n-r \\
0 & \mbox{if } j=n-r
\end{array}\right.
$$
and
$$
J/I \cong S_{\frac{r-t}{2}} (M_{\ffi}) \otimes S_{\frac{r-t}{2}}(P) \otimes
\wedge^f F^* \otimes \wedge^g G \otimes
\wedge^t P.
$$
\end{itemize}
\end{proposition}

\begin{proof} For $i = 0,\ldots,r$ the module $S_i(M_{\ffi})$ is a perfect
  module of dimension $n-r$. Hence claim (a) follows by
Proposition~\ref{cohomology-of locus} and Lemma~\ref{unmixedness-criterion}.

In order to show (b) we note first that the maximal codimension of a
component of I is $r+1$ because $\depth R/I = n-r$. Let $Q$ be the
intersection of all these components and let $J$ be the intersection of the
remaining ones. Then $I=J\cap Q$ and the components of $J$
have codimension $\leq r$. We have to show that $J$ is
unmixed.

As a first step we will prove that $\depth R/J > n-r$. We induct on $n-r\geq
0$. If $r=n$ the claim is clear since $J$ is saturated by construction.

Let $r=n-1$. It follows that  $\dim R/I = 1 + t \geq2$ and that $R/Q$ and
$S_{\frac{r+t}{2}}(\Mf)$ have positive dimension.
Now we look at the exact sequence
$$
0 \to R/I \to R/J \oplus R/Q \to R/ (J + Q) \to 0.
$$
We can find an $l \in [R]_1$ which is a parameter on $R/I, R/J, R/Q,
S_{\frac{r+t}{2}}(\Mf)$ and also on $R/(J+Q)$ if it has positive
dimension. Using $\HH^2(R/I) \cong \HH^2(R/J)$ we obtain a commutative
diagram with exact rows
$$
\begin{array}{cccccccc}
\HH^1(R/I)(-1)  & \to & \HH^1(R/J)(-1) & \oplus & \HH^1(R/Q)(-1) & \to &
\HH^1(R/(J+Q))(-1) & \to  0 \\
\downarrow\!{\scriptstyle \beta_1} & & \downarrow\!{\scriptstyle \beta_2} &
& \downarrow\!{\scriptstyle \beta_3} & & \downarrow\!{\scriptstyle \beta_4}
& \\
\HH^1(R/I) & \to & \HH^1(R/J) & \oplus & \HH^1(R/Q) & \to & \HH^1(R/J+Q) & \to
0
\end{array}
$$
where the vertical maps are multiplication by $l$.  Due to our choice of
$l$ it holds $\HH^1(R/(Q + l R)) = \HH^1(R/(J+Q+l R)) = 0$. Thus $\beta_3$ and
$\beta_4$ are surjective. Since $l$ is a parameter of the Cohen-Macaulay
module $S_{\frac{r+t}{2}}(M)$ the multiplication map
$S_{\frac{r+t}{2}}(M)(-1)
\stackrel{l}{\longrightarrow} S_{\frac{r+t}{2}}(M)$ is injective, thus the
dual map $S_{\frac{r+t}{2}}(M)^{\vee}
\stackrel{l}{\longrightarrow} S_{\frac{r+t}{2}}(M)^{\vee}(1)$ is an
epimorphism. Therefore $\beta_1$ is surjective due to
Proposition~\ref{cohomology-of locus}. The same is true for $\beta_2$ by the
commutative diagram above. Since $R/J$ is unmixed
Lemma~\ref{unmixedness-criterion} implies that $\HH^1(R/J)$ is finitely
generated, hence it must be zero by Nakayama's lemma.

Finally, let $r \leq n-2$. We consider the commutative diagram
$$
\begin{array}{ccccccccc}
0 \to & R/I(-1)  & \to & R/J(-1) & \oplus & R/Q(-1) & \to & R/(J+Q)(-1) & \to
0 \\
& \downarrow & & \downarrow & & \downarrow\ & & \downarrow  & \\
0 \to & R/I & \to & R/J & \oplus & R/Q & \to & R/(J+Q) & \to
0
\end{array}
$$
where the vertical maps are multiplication by $l$. By Lemma~\ref{depth-formula}
$\depth R/I \geq 2$ and by assumption on $r$ $R/J$ and
$R/Q$ have positive depth. Hence the cohomology sequence induced by the
bottom line provides $\HH^0(R/(J+Q)) = 0$. It follows that we may choose $l$ as
non-zero divisor on $R/I, R/J, R/Q$ and $R/(J+Q)$, i.e., the vertical maps in
the diagram above are all injective. Thus the Snake lemma implies the exact
sequence
$$
0 \to \overline{R}/\overline{I} \to \overline{R}/\overline{J} \oplus
\overline{R}/\overline{Q} \to \overline{R}/\overline{J} + \overline{Q} \to
0
$$
where we denote by $^-$ again the functor $\_\_ \otimes_{R} R/l R$. By
Bertini's theorem $\overline{Q}$ is unmixed of codimension $r+1$ in
$\overline{R}$ (possibly) up to a component associated to the irrelevant
ideal of  $\overline{R}$. Moreover, we have seen in the proof of
Lemma~\ref{depth-estimate} that the induction hypothesis applies to
$\overline{I}$. But the last exact sequence implies $\overline{I} =
\overline{J} \cap \overline{Q}$. It follows that $\overline{J}$ is the
intersection of the top-dimensional components of $\overline{I}$. Hence we
get by induction $\depth \overline{R}/\overline{J} \geq n-r$. But $l$ was a
non-zero divisor on $R/J$ thus we obtain
$$
\depth R/J > \depth \overline{R/J} = \depth \overline{R}/\overline{J}
$$
completing our induction.

Next we consider the exact sequence
$$
0 \to J/I \to R/I \to R/J \to 0. \leqno(*)
$$
Note that $J/I\cong (J+Q)/Q$ has dimension $\leq n-r$. Moreover we have just
shown $\depth R/J > n-r$. Therefore the induced cohomology sequence yields:
$$
\HH^i(R/J) \cong \HH^i(R/I) \quad \mif i > n-r,
$$
$$
\HH^{n-r}(J/I)   \cong  \HH^{n-r} (R/I)  \cong S_{\frac{r+t}{2}}
(M_{\ffi})^\vee\otimes S_{\frac{r-t}{2}} (P) \otimes \wedge^t P (1-r(R))
$$
and $J/I$ is Cohen-Macaulay of dimension $n-r$. The first
isomorphisms, Proposition~\ref{cohomology-of locus} and
Lemma~\ref{unmixedness-criterion} imply now that $J$ must be unmixed.

Now we use that for a Cohen-Macaulay module $M$ it holds $K_{K_M} \cong
M$. Thus we get by duality and Lemma~\ref{can-module-of-symm-powers}
\begin{eqnarray*}
J/I & \cong & \HH^{n-r}( \HH^{n-r}(J/I)^{\vee})^{\vee} \\
&  \cong & \HH^{n-r}(S_{\frac{r+t}{2}} (M_{\ffi}))^\vee \otimes
S_{\frac{r-t}{2}} (P)  \otimes \wedge^t P (1 - r(R)) \\
& \cong & \Ext^{r+1}(S_{\frac{r+t}{2}} (M_{\ffi}),R) \otimes
S_{\frac{r-t}{2}} (P)  \otimes \wedge^t P \\
& \cong & S_{\frac{r-t}{2}} (M_{\ffi}) \otimes S_{\frac{r-t}{2}} (P) \otimes
\wedge^f F^* \otimes \wedge^g G.
\otimes \wedge^t P.
\end{eqnarray*}
This finishes the proof.
\end{proof}

\begin{remark} 
The arguments in the previous proof also provide that $R/Q$ is
Cohen-Macaulay of dimension $n-r$ and
$$
n-r-1 \leq \depth R/(J+Q) \leq \dim R/(J+Q) \leq n-r.
$$
\end{remark}

Our results with respect to ring-theoretic properties can be summarized as
follows. Recall that $X$ denotes the top-dimensional part of the degeneracy
locus $S$ of $\psi$. 

\begin{theorem} \label{summary_for_degeneracy_loci}
With the  notation above we have:
\begin{itemize}
\item[(a)] If $r = n$ then $S$ is equidimensional and locally
  Cohen-Macaulay.
\item[(b)] $S$ is equidimensional if and only if $r+t$
  is odd or $r = n$.

  Moreover, if $r < n$ then $X$ is locally Cohen-Macaulay if and only if $X$ is
  arithmetically Cohen-Macaulay.
\item[(c)] If $r+t$ is odd then $X = S$ is arithmetically Cohen-Macaulay if
  and only if $t=1$. In this case $S$ has Cohen-Macaulay type $\leq 1 +
  \binom{\frac{r}{2}+g-1}{g-1}$.
\item[(d)] Let  $r+t$ be  even. Then
\begin{itemize}
\item[(i)] $X$ is arithmetically Cohen-Macaulay if and only if $1 \leq t
\leq 2$. If $t =1$ then $X$ is arithmetically Gorenstein. If $t = 2$ then $X$
has Cohen-Macaulay type $\leq r - 1 + \binom{\frac{r}{2}+g-1}{g-1} \cdot
(\frac{r}{2} - 1)$.
\item[(ii)] If in addition  $r < n$ then the components
  of $S$ have either codimension $r-t+1$ or codimension $r+1$.
\end{itemize}
\end{itemize}
\end{theorem}

\begin{proof} (a) If $r = n$ then $\dim \Mf = 0$. Hence it follows by
  Proposition~\ref{cohomology-of locus} that the modules $H^i_*(Z,\cJ_S)$
  have finite length if $i \leq \dim S$ which is equivalent to $S$ being
  equidimensional and locally Cohen-Macaulay.

(b)  If $r+t$ is odd then $S$ is equidimensional due to
Proposition~\ref{coho_of_top-dimensional}. If $r+t$ is even then
Proposition~\ref{coho_of_top-dimensional} shows that $S$ is 
equidimensional if and only if $r=n$.

Moreover, $r<n$ implies that $\Mf$ does not have finite length. Therefore
Proposition~\ref{coho_of_top-dimensional} furnishes that $R/J(\psi)$ has
cohomology of finite length if and only if $R/J(\psi)$ is
Cohen-Macaulay. Claim (b) follows.

(c) If $r+t$ is odd then we have by Proposition~\ref{coho_of_top-dimensional}
that $S$ is equidimensional and by Corollary~\ref{depth-formula} that $\depth
R/I(\psi) = n-r+1$. The claim follows because of $\dim R/I(\psi) = n-r+t$.

(d) If $r+t$ is even we get $r \leq t-2$. Hence the claim is a
consequence of Propositions~\ref{cohomology-of locus} and
\ref{coho_of_top-dimensional}.
\end{proof}

\begin{remark} (i) If we specialize the previous result to $Z = \PP^n, r=3,
  t=1$ then we get the main result of \cite{Mig-P_gorenstein}. \\
(ii) In case $t=1$ and $r$ is even the result has been first proved by the
first and third author who communicated it to Kustin. Subsequently, Kustin
\cite{Kustin} strengthened it by removing almost all the assumptions on the
ring $R$ and computing a free resolution (cf.\ also Remark
\ref{discussion_of_minimality} and Corollary \ref{one-sect-even-case}).  \\ 
(iii) Note that in the above theorem the term equidimensional is used in the 
scheme-theoretic sense. Thus $S$ being equidimensional  does not automatically
imply that $I(\psi)$ is saturated. However, 
due to Corollary 4.3 $I(\psi)$ is not saturated if and 
only if $r = n$ and $r+t$ is even. We use this fact in
Section~\ref{resolution_general}. 
\end{remark}

The next result generalizes \cite{Mig-P_gorenstein}, Corollary 1.2.

\begin{corollary} Let $\cE$ be a vector bundle on $\PP^n$ of rank $n$ where
  $n$ is odd such that $H^i_*(\PP^n,\cE) = 0$ if $2 \leq i \leq n-1$.
Let $s$ be a section of $\cE$ vanishing on a scheme $X$ of codimension
$n$. Then $X$ is  arithmetically Gorenstein.
\end{corollary}

\begin{proof} According to Corollary~\ref{Bu-Rim-sheafs-on-proj-space},
  $\cE$ is a Buchsbaum-Rim sheaf. Therefore
Theorem~\ref{summary_for_degeneracy_loci} shows the claim.
\end{proof}


\section{The resolutions of the loci}
\label{resolution_general} 

In this section  we show how free resolutions can be obtained
for the schemes described in this paper. In most cases we expect these
resolutions to be minimal. The main tools are the Eagon-Northcott complex
$E_{\bullet}$ in $(4.1)$, 
its dual and a general result for comparing the resolution of the scheme
with its cohomology modules. 

We begin by considering the Eagon-Northcott complex again.  The interested
reader will have observed that this complex is {\it not} exact in
general. Fortunately, we are able to compute its homology.

\begin{proposition} \label{homology_of_EN} The homology modules of the
  Eagon-Northcott $E_{\bullet}$ complex are:
$$
H_i (E_{\bullet}) \cong \left\{ \begin{array}{ll}
S_j(\Mf) \otimes S_{r-t-j}(P) \otimes \wedge^f F^* \otimes \wedge^g G &
\mif  i=r-1-2j
\mbox{ where } j \in \ZZ, \; t \leq i \leq r-3 \\
0 & \mbox{otherwise}.
\end{array}\right.
$$
\end{proposition}

\begin{proof} According to Lemma~\ref{prop-of-exteriour-powers} and
Lemma~\ref{depth-estimate} all non-trivial modules occurring in
  $E_{\bullet}$ have depth $\geq n-r+1$. Thus using $n-r$ general
  linear forms in $R$ and arguments as in Proposition~\ref{depth-estimate}
  we see that it suffices to consider the case where $R$ has dimension
  $r+1$, i.e., we may and will assume that $n=r$.

Then we already know that the
homology modules of $E_{\bullet}$ have finite length. As for the
cohomology of $R/I$ one computes (cf.\ Proposition~\ref{cohomology-of
  locus})
$$
H^1_m(\im \delta_i) \cong \left\{ \begin{array}{ll}
S_{r-j} (\Mf)^\vee \otimes S_{r-j-t}(P)(1- r(R))  & \mif t \leq i = r-2j <
r \\
0 & \mif  r+i \mbox{ is odd}
\end{array} \right.
$$
Since $r = n$ the module $\Mf$ has finite length. It follows
\begin{eqnarray*}
S_{r-j}(\Mf)^{\vee}(1-r(R)) & \cong & \HH^0(S_{r-j}(\Mf))^{\vee}(1-r(R)) \\
& \cong & \ER^{r+1}(S_{r-j}(\Mf),R) \quad \mbox{by \ref{duality}} \\
& \cong & S_j(\Mf) \otimes \wedge^f F^* \otimes \wedge^g G \quad \mbox{by
  \ref{can-module-of-symm-powers}.}
\end{eqnarray*}
Thus we have
$$
H^1_m(\im \delta_i) \cong \left\{ \begin{array}{ll}
S_{j} (\Mf) \otimes S_{r-j-t}(P) \otimes \wedge^f F^* \otimes \wedge^g
G  & \mif t \leq i = r-2j <  r \\
0 & \mif  r+i \mbox{ is odd}
\end{array} \right. \leqno(+)
$$
The modules $\im \delta_i$ and $\ker \delta_i$ are submodules
of the reflexive modules $E_{i+1}$ and $E_i$, respectively.
Hence both have positive depth.
Since $\depth E_i = n+1-i \geq 2$ if $i<r=n$ the exact
sequence
$$
0 \to \ker \delta_i \to E_i \to \im \delta_i \to 0
$$
shows $\depth \ker \delta_i \geq 2$ for all $i=t, \ldots, r$.

Using the finite length of $H_i(E_{\bullet})$ the exact sequence
$$
0 \to \im \delta_{i+1} \to \ker \delta_i \to H_i (E_{\bullet} ) \to 0
$$
provides $H^1_m (\im \delta_{i+1}) \cong H^0_m (H_i(E_{\bullet}))
\cong H_i(E_{\bullet})$. Hence $(+)$ proves our assertion because
$\delta_t$ is surjective by the definition of the ideal $I$. 
\end{proof}

\begin{remark} \label{remark_resolution_in_general} Now we can compute a
  free resolution of $I = I(\psi)$ as follows: Proposition~\ref{homology_of_EN}
provides exact sequences
$$
0 \to \ker \delta_{r+1-2j} \to E_{r+1-2j} \to E_{r-2j} \to
\im \delta_{r-2j} \to 0, \quad t \leq r-2j<r, \leqno(1)
$$
and if $r+t$ is odd additionally
$$
0 \to \ker \delta_t \to E_t \to I(p) \to 0,
\leqno(1')
$$
$$
0 \to \im \delta_{r-2j} \to \ker \delta_{r-2j-1} \to S_j(\Mf) \otimes
S_{r-t-j}(P) \otimes \wedge^f F^* \otimes \wedge^g G \to 0, \quad
t<r-2j<r. \leqno(2)
$$
According to Proposition~\ref{EN-complexes_are_exact} and
Proposition~\ref{prop-of-exteriour-powers} we know the minimal
resolution of
$E_i$ and $S_j(\Mf)$, respectively. Thus, in order to get a resolution of
$I$ we
compute  successively resolutions of $\im \delta_{r-2},  \ker
\delta_{r-3}, \\
\im \delta_{r-4}, \ldots, \im \delta_t = I(p)$
using the exact sequences (1) and (2). If (1) is used we just
apply the mapping cone procedure twice. If we use (2) we apply the
Horseshoe lemma. 
\end{remark}

Following the procedure just described  we certainly do not obtain a
minimal free 
resolution of $I$. Thus we need some results which allow us to split off
redundant terms. The idea is to compare the resolution of $I$ with those of
its cohomology modules. In particular,  this requires information on the
canonical module $K_S = \ER^{r-t}(R/I(\psi),R) (r(R) - 1)$ of our
degeneracy locus $S$ which we derive first.  

\begin{lemma} \label{can_module_depth_estimate} The depth of the canonical
  module of 
  $X$ is at least $\min \{n-r+t, n-r+2\}$. In particular, it is
  Cohen-Macaulay if $1 \leq t \leq 2$.  
\end{lemma} 

\begin{proof} Let $A = R/I(\psi)$. We induct on $n - r$. Since
  the canonical module always satisfies $\depth K_A \geq \min \{\dim A,
  2\}$ the claim is clear if $n = r$. 

Let $n > r$ and let $l \in R$ be a general linear form. We want to show 
$$
K_A/ l K_A \cong K_{A/l A}. 
$$
Due to Corollary \ref{depth-formula} there is an exact sequence induced by
multiplication 
$$
0 \to A(-1) \stackrel{l}{\longrightarrow} A \to A/l A \to 0. 
$$
It provides the long exact sequence 
$$
0 \to \ER^{r-t}(A,R) \stackrel{l}{\longrightarrow} \ER^{r-t}(A,R)(1) \to
\ER^{r-t+1}(A/l A,R) \to \ER^{r-t+1}(A,R) \stackrel{l}{\longrightarrow}
\ER^{r-t+1}(A,R)(1) \to \ldots. 
$$
Using $\ER^{r-t-1}(A/l A,R)(r(R) - 1)  \cong K_{A/l A}$ we can rewrite the
last sequence as 
$$
0 \to K_A(-1) \stackrel{l}{\longrightarrow} K_A \to K_{A/l A} \to
\ER^{r-t+1}(A,R)(r(R)-2) \stackrel{l}{\longrightarrow} 
\ER^{r-t+1}(A,R)(r(R)-1) \to \ldots. 
$$
We claim that the multiplication map on the right-hand side is
injective. Indeed, if $r+1$ is odd or $2t \geq r+1$ then we have by duality
and Proposition 
\ref{cohomology-of locus} that $\ER^{r-t+1}(A,R) = 0$. Otherwise  the
multiplication map is (up to degree shift) essentially given by 
$$
S_{\frac{r+1}{2}}(\Mf) \otimes S_{\frac{r+1}{2}-t} (P)^* (-1)
\stackrel{l}{\longrightarrow} S_{\frac{r+1}{2}}(\Mf) \otimes
S_{\frac{r+1}{2}-t} (P)^*. 
$$
But the module $S_{\frac{r+1}{2}}(\Mf)$ is Cohen-Macaulay of dimension $n-r
> 0$ according to Lemma \ref{prop-of-exteriour-powers}. Thus we can choose
$l$ as a regular element with respect to $S_{\frac{r+1}{2}}(\Mf)$ and the
claim follows. 

Hence the exact sequence above implies 
$$
K_A/ l K_A \cong  K_{A/l A} \quad \mbox{and} \quad \depth K_A > \depth
K_{A/l A}. 
$$
Applying the induction hypothesis to $K_{A/l A}$ completes the proof. 
\end{proof} 

\begin{remark} \label{X_and_S_have_same_can_module} The canonical modules
  of $S$ and its top-dimensional part $X$ are isomorphic. This follows from
  the corresponding (slightly stronger) result for $R/I(\psi)$ and
  $R/J(\psi)$. Indeed, according to Proposition
  \ref{coho_of_top-dimensional} there is an exact sequence 
$$
0 \to S_{\frac{r-t}{2}}(\Mf) \otimes S_{\frac{r-t}{2}}(P) \otimes \wedge^f
F^* \otimes \wedge^g G \otimes \wedge^t P \to R/I(\psi) \to R/J(\psi) \to
0. 
$$
Since $\dim S_{\frac{r-t}{2}}(\Mf) = n-r < n-r+t = \dim R/I(\psi) = \dim
R/J(\psi)$ the claim follows by the long exact cohomology sequence. 
\end{remark} 

\begin{lemma} \label{dual_of_EN_complex} The dual of the Eagon-Northcott
  complex $E_{\bullet}$ provides a complex 
$$
0 \to \wedge^t P \to E_t^* \stackrel{\delta^*_{t+1}}{\longrightarrow}  \ldots \to
  E_{r-1}^* \stackrel{\delta^*_{r}}{\longrightarrow} E_r^*
  \stackrel{\gamma}{\longrightarrow} K_{R/I} \otimes \wedge^t P (1 - r(R))
  \to 0 
$$ 
which we denote (by slight abuse of notation) by $E_{\bullet}^*$. Its
(co)homology modules are given by 
$$
H^i(E_{\bullet}^*) \cong \left \{ 
\begin{array}{ll} 
S_j(\Mf) \otimes S_{j-t}(P)^* & \mif 2t + 1 \leq i = 2 j + 1 \leq r+1 \\
0 & \mbox{otherwise} 
\end{array} 
\right. 
$$ 
In particular, $E_{\bullet}^*$ is exact if $t \geq \frac{r+1}{2}$. 
\end{lemma} 

\begin{proof} If $r = 2$ (and thus $t=1$) the claim follows immediately by 
  dualizing because $E_{\bullet}$ is exact in this case. 

Now let $r \geq 3$. We begin by verifying that $E_{\bullet}^*$ is indeed a
complex. This is only an issue in the beginning of this sequence. 
Consider the exact sequence 
$$
0 \to E_r \stackrel{\delta_r}{\longrightarrow} E_{r-1} \to \im \delta_{r-1}
\to 0. 
$$
Dualizing provides the exact sequence 
$$
 E_{r-1}^* \stackrel{\delta^*_{r}}{\longrightarrow} E_r^*
\stackrel{\alpha}{\longrightarrow} \ER^1(\im \delta_{r-1}, R) \to
\ER^1(E_{r-1},R). \leqno(+) 
$$
The module on the right-hand side vanishes since $\HH^n(\wedge^{r-1}\Bf^*)
= 0$ if $r \neq 2$. Hence the map $\gamma$ is surjective. In order to
compare its image with $K_{R/I}$ we look at the proof of Proposition
\ref{cohomology-of locus} and use its notation. 

If $t \geq \frac{r+1}{2}$ we have (cf.\ Case 1) 
$$
\HH^{n+t-r}(R/I)(p) \cong \HH^n(\im \delta_{r-1}). 
$$
Thus $(+)$ provides the exact  sequence 
$$
 E_{r-1}^* \stackrel{\delta^*_{r}}{\longrightarrow} E_r^*
\stackrel{\alpha}{\longrightarrow} K_{R/I}(1 - r(R) -p) \to 0.
$$
Putting $\gamma = \alpha$ we get the desired complex and 
$$
H^{r+1}(E_{\bullet}^*) = H^{r}(E_{\bullet}^*) = 0. 
$$
Now let $t \leq \frac{r+1}{2}$ and let $r$ be odd. Then we have (cf.\ Case
2) an embedding 
$$ 
\beta: \HH^{n+t-r}(R/I)(p) \hookrightarrow  \HH^n(\im \delta_{r-1}).
$$
We define $\gamma$ as the composition of $\alpha$ and $\beta^{\vee}$ (with
the appropriate shift). Thus $\gamma$ is surjective, i.e. 
$$
H^{r+1}(E_{\bullet}^*) =  0. 
$$ 
Finally, let $t \leq \frac{r+1}{2}$ and let $r$ be even. Then we have (cf.\
Case 3) an exact sequence 
$$
0 \to \HH^{n+1-\frac{r}{2}}(E_{\frac{r}{2}}) \to \HH^{n+1-\frac{r}{2}}(\im
\delta_{\frac{r}{2}}) \to  \HH^{n+2-\frac{r}{2}}(\ker
\delta_{\frac{r}{2}}) \to \HH^{n+2-\frac{r}{2}}(E_{\frac{r}{2}}) = 0 
$$
and isomorphisms
$$
\HH^{n+1-\frac{r}{2}}(\im \delta_{\frac{r}{2}}) \cong \HH^{n-r+t}(R/I)(p),
\quad 
\HH^{n+2-\frac{r}{2}}(\ker \delta_{\frac{r}{2}}) \cong \HH^n(\im
\delta_{r-1}). 
$$
Thus we can conclude with the help of $(+)$ that there is an exact sequence
$$
E_{r-1}^* \stackrel{\delta^*_{r}}{\longrightarrow} E_r^*
  \stackrel{\gamma}{\longrightarrow} K_{R/I} (1 - p - r(R)) \to
  \HH^{n+1-\frac{r}{2}}(E_{\frac{r}{2}})^{\vee}(1-r(R)) 
  \to 0. 
$$
Using Lemma \ref{prop-of-exteriour-powers} it follows 
$$
H^{r+1}(E_{\bullet}^*) \cong
\HH^{n+1-\frac{r}{2}}(E_{\frac{r}{2}})^{\vee}(1-r(R)) \cong
S_{\frac{r}{2}}(\Mf) \otimes S_{\frac{r}{2}-t}(P)^*. 
$$
Hence we have shown in all cases that $E_{\bullet}^*$ is a complex and we
have computed $H^{r+1}(E_{\bullet}^*)$. 

In order to compute the other cohomology modules we proceed as in the
proofs of Proposition \ref{cohomology-of locus} and Proposition
\ref{homology_of_EN}. Indeed, we just have to  use Lemma
\ref{can_module_depth_estimate} as replacement of Lemma
\ref{depth-estimate} and $E_{\bullet}^*$ instead of $E_{\bullet}$. The
details are tedious but straightforward. We omit them. 
\end{proof} 

The next result is also interesting in its own right. It relates the minimal free resolution of a module
to those of (the duals) of its cohomology modules. Observe that the result 
as it is stated remains valid even if $R$ is not a Gorenstein but just a
Cohen-Macaulay ring, though we won't use this fact here. 

\begin{proposition} \label{resolution_via_cohomology} Let $N$ be a finitely
  generated graded torsion $R$-module which has projective dimension
  $s$. Then it holds for all integers $j \geq 0$  that
  $\TR_{s-j}(N,K)^{\vee}$ is a direct summand of 
$$
\oplus_{i=0}^j \TR_{j-i}(\ER^{s-i}(N,R),K). 
$$ 
Moreover, we have $\TR_{s}(N,K)^{\vee} \cong \TR_{0}(\ER^{s}(N,R),K)$ and
that $\TR_{1}(\ER^{s}(N,R),K)$ is a direct summand of
$\TR_{s-1}(N,K)^{\vee}$. 
\end{proposition} 

\begin{proof} For the purpose of this proof we write $M \ds N$ in order to
  express that the submodule $M$ is a
  direct summand of the module $N$. 

Consider a minimal free resolution of $N$: 
$$
0 \to F_s \to \ldots \to F_1 \to F_0 \to N \to 0. 
$$ 
Dualizing with respect to $R$ provides the complex: 
$$
0 \to F_0^* \stackrel{\alpha_0}{\longrightarrow} F_1^* \to \ldots \to
F_{s-1}^* \stackrel{\alpha_{s-1}}{\longrightarrow} F_s^* \to \ER^s(N.R) \to
0. 
$$
Since the maps $\alpha_j$ are duals of minimal maps we can write 
$$
\ker \alpha_j \cong G_j \oplus M_j 
$$ 
where $M_j$ does not have a free $R$-module as direct summand  and
$G_j$ is a free $R$-module being a direct summand of $F_j^*$ (but possibly
trivial).  
Moreover, there are exact sequences 
$$
0 \to \im \alpha_{j-1} \to \ker \alpha_j \to \ER^j(N,R) \to 0. 
$$
Since $\alpha_{j-1}$ is a minimal homomorphism it holds $\im \alpha_{j-1} \subset \fm
\cdot F_j^*$. This shows: \\
(1) The minimal generators of $G_j$ give rise to minimal generators of
$\ER^j(N,R)$. 

Now we consider the diagram 
$$
\begin{array}{ccccccc}
0 \to & \im \alpha_{j-1} & \to & \ker \alpha_j & \to & \ER^j(N,R) & \to 0\\
& \downarrow & & \downarrow & & & \\
& F^*_j & = & F^*_j & & &
\end{array}
$$ 
where the vertical maps are the natural embeddings. Thus it is commutative
and  we obtain an exact sequence 
$$
0 \to \ER^j(N,R) \to F_j^*/\im \alpha_{j-1} \to F_j^*/ \ker \alpha_j \to
0. 
$$
The Horseshoe lemma yields a free resolution of the middle module as direct
sum of the resolutions of the outer modules. After splitting off redundant
terms we get a minimal free resolution, i.e, using $\TR_i(\im
\alpha_{j-1},K) \cong \TR_{i+1}(F_j^*/\im \alpha_{j-1},K)$  
we obtain  
$$
\TR_i(\im \alpha_{j-1},K) \ds \TR_{i+1}(\ER^j(N,R),K) \oplus \TR_i(\ker
\alpha_j,K) \quad (i \geq 0). \leqno(2)
$$
Now we are ready to show by induction on $j \geq 1$: 
$$
\TR_i(\im \alpha_{s-j}, K) \ds \oplus_{k=0}^{j-1}
\TR_{i+j-k}(\ER^{s-k}(N,R), K) \quad (i \geq 0) \leqno(*)
$$
and 
$$
\TR_{s-j+1}(N, K)^{\vee} \ds \oplus_{k=0}^{j-1}
\TR_{j-1-k}(\ER^{s-k}(N,R), K). \leqno(**)
$$ 
Let $j = 1$. Since $\alpha_{s-1}$ is a minimal homomorphism the exact
sequence 
$$
0 \to \im \alpha_{s-1} \to F_s^* \to \ER^s(N,R) \to 0 
$$
implies 
$$ 
\TR_{s}(N,K)^{\vee} \cong \TR_{0}(\ER^{s}(N,R),K),  
$$ 
thus in particular $(**)$, and 
$$
\TR_{i}(\im \alpha_{s-1}, K) \cong \TR_{i+1}( \ER^s(N,R), K) \quad (i \geq
0) 
$$
which shows $(*)$. 

Let $j \geq 2$. Consider the exact sequence 
$$
0 \to \ker \alpha_{s-j+1} \to F_{s-j+1}^* \to \im \alpha_{s-j+1} \to 0.
\leqno(3) 
$$
We obtain by the definition of $G_{s-j+1}$ 
\begin{eqnarray*}
F_{s-j+1}^* \otimes K & \cong & (G_{s-j+1} \otimes K) \oplus \TR_0(\im
\alpha_{s-j+1}, K) \\[3pt]
& \ds & \TR_0(\ER^{s-j+1}(N,R), K) \oplus \bigoplus_{k=0}^{j-2}
\TR_{j-1-k}(\ER^{s-k}(N,R), K) \quad \mbox{(by (1) and induction)} \\[3pt]
& = &   \bigoplus_{k=0}^{j-1} \TR_{j-1-k}(\ER^{s-k}(N,R), K) 
\end{eqnarray*} 
which shows $(**)$. 

Furthermore, (3) provides  
\begin{eqnarray*} 
\TR_i(\ker \alpha_{s-j+1}, K) & \ds & \TR_{i+1}(\im \alpha_{s-j+1}, K) \\[3pt]
& \ds & \bigoplus_{k=0}^{j-2} \TR_{i+j-k}(\ER^{s-k}(N,R), K) \quad \mbox{(by
  induction).} 
\end{eqnarray*}
Thus we conclude with the help of (2):
\begin{eqnarray*} 
\TR_i(\im \alpha_{s-j}, K) & \ds & \TR_{i+1}(\ER^{s-j+1}(N,R), K) \oplus
\bigoplus_{k=0}^{j-2} \TR_{i+j-k}(\ER^{s-k}(N,R), K) \\[3pt]
& = & \bigoplus_{k=0}^{j-1} \TR_{i+j-k}(\ER^{s-k}(N,R), K). 
\end{eqnarray*}
This proves $(*)$, i.e., the induction step is complete. 

It only remains to verify the very last assertion. Looking at (3) with
$j=2$ we get
$$
\TR_0(\im \alpha_{s-1}, K) \ds F_{s-1}^* \otimes K. 
$$
But at the beginning of the induction we have seen that $\TR_0(\im
\alpha_{s-1}, K) \cong \TR_1(\ER^s(N,R), K)$. Our claim follows 
\end{proof} 

\begin{remark} \label{comparison_with_Rao} In order to explain how the last
  result can be used we compare it with a well-known statement of Rao. For
  this let $R = K[x_0,\ldots,x_3]$ and let $I_C \subset R$ be the
  homogeneous ideal of a curve $C \subset \PP^3$. Put $A = R/I_C$ and
  consider the following  minimal free resolutions:
\[
\begin{array}{c}
0 \rightarrow F_3  \rightarrow F_2 \rightarrow F_1 \rightarrow R
\rightarrow A \rightarrow 0 \\ \\
0 \rightarrow G_4 \rightarrow G_3 \rightarrow G_2   \rightarrow G_1
\rightarrow G_0 \rightarrow \ER^3(A,R)  \rightarrow 0 \\ \\
0 \rightarrow D_2 \rightarrow D_{1} \rightarrow D_{0} \rightarrow
\ER^2(A,R) \rightarrow 0. 
\end{array}
\] 
Then the additions to our general observation in the lemma above yield 
$$
G_0 \cong F_3^* \quad \mbox{and} \quad G_1 \ds F_2^*. 
$$ 
This is precisely the content of Rao's Theorem 2.5 in \cite{R1}. Our Lemma
\ref{resolution_via_cohomology} gives in addition  
$$
G_1 \ds F_2^* \ds G_1 \oplus D_0 
$$
and 
$$ 
F_1^* \ds G_2 \oplus D_1 
$$  
because $\ER^1(A,R) = 0$. 
\end{remark} 

Now we have all the tools for establishing the main result of this section. 

\begin{theorem} \label{resolution_of_top-dim_part} Consider the following
  modules where we use the conventions that $i$ and $j$ are non-negative
  integers and that a sum is trivial if it has no summand: 
$$
A_k = \bigoplus_{\begin{array}{c} 
{\scriptstyle i+2j = k + t -1}\\ [-4pt]
{\scriptstyle t \leq i+j \leq \frac{r+t-1}{2}}
\end{array}} 
\wedge^i F^* \otimes S_j(G)^* \otimes S_{i+j-t}(P),  
$$
$$
C_k = \bigoplus_{\begin{array}{c} 
{\scriptstyle i+2j = r+1-t-k}\\ [-4pt]
{\scriptstyle i+j \leq \frac{r-t}{2}}
\end{array}} 
\wedge^i F \otimes S_j(G) \otimes S_{r-t-i-j}(P) \otimes \wedge^f F^*
\otimes \wedge^g G.  
$$ 
Observe that it holds: 

$A_r = 0$ if and only if $r+t$ is even, 

$C_1 = 0$ if and only if $r+t$ is odd,

$C_k = 0$ if $k \geq r+2-t$  and 

$C_{r+1-t} = S_{r-t}(P) \otimes \wedge^f F^* \otimes \wedge^g G$. \\
Then the homogeneous ideal $I_X = J(\psi)$ of the top-dimensional part $X$
of the degeneracy locus $S$ has a graded free resolution of the form 
$$
0 \to A_r \oplus C_r \to \ldots \to A_1 \oplus C_1 \to I_X \otimes \wedge^t
P^* \to 0. 
$$
\end{theorem}

\begin{proof} Following the procedure described in Remark
  \ref{remark_resolution_in_general} we get a  resolution of
  $I(\psi)$. According to Proposition \ref{coho_of_top-dimensional} it
  holds 
  $I(\psi) = J(\psi)$ if and only if $r+t$ is odd. If $r+t$ is even we have
  by the same result an exact sequence 
$$
0 \to I(\psi) \to J(\psi) \to S_{\frac{r-t}{2}}(\Mf) \otimes
S_{\frac{r-t}{2}}(P) \otimes \wedge^f F^* \otimes \wedge^g G \otimes
\wedge^t P^*.  
$$
Thus, using the Horseshoe lemma we get a finite free resolution of
$J(\psi)$ in any case. It is not minimal. In order to split  off redundant
terms we proceed as follows: 

Since the resolution is finite the Auslander-Buchsbaum formula and
Proposition  \ref{coho_of_top-dimensional} yield the projective dimension
of $J(\psi)$. Thus, in a first step we can split off all the terms in the
resolution of $J(\psi)$ occurring past its projective dimension. 

Next, we use Lemma \ref{dual_of_EN_complex} (cf.\ Remark
\ref{X_and_S_have_same_can_module}) in order to obtain a free resolution of
$\ER^{r-t+1}(R/J(\psi),R)$. For $j \neq r-t+1$ we know a free
resolution of $\ER^{j}(R/J(\psi),R)$ by Proposition
\ref{coho_of_top-dimensional} and Proposition
\ref{EN-complexes_are_exact}. Hence, in a second step we can split off
further terms in the resolution of $J(\psi)$  by applying Proposition
\ref{resolution_via_cohomology}. This provides the resolution as
claimed. The details are very tedious but straightforward. We omit them. 
\end{proof} 

Since the above proof is somewhat sketchy, we illustrate it by deriving the
resolution for $t = 2$ and $r = 6,7$ (cf.\ also the next corollaries). 

\begin{example} \label{t=2_r=6,7} Using our standard notation we define
  integers $c, p$ by 
$$
R(c) \cong \wedge^f F \otimes \wedge^g G^* \quad \mbox{and} \quad R(p) =
\wedge^t P^*. 
$$
(i) Let $t=2$ and $r=7$. Then we know that $I = I(\psi)$ is unmixed and
that $X$ is not arithmetically Cohen-Macaulay. Using the Eagon-Northcott
complex $E_{\bullet}$ we see that a free resolution of $I$ begins as
follows: 
$$
\begin{array}{c@{\ }c@{\ }c@{\ }c@{\ }c@{\ }c@{\ }c@{\ }c@{\ }c@{\ }c@{\
      }c@{\ }c@{\ }c@{\ }c@{\ }c} 
&&&&& \downarrow \\
\wedge^2 F^* \otimes S_2(G)^* \otimes S_2(P) & \oplus & 
 S_3(G)^* \otimes P  &  & & \oplus & & & 
\wedge^3 F^* \otimes S_2(P) \\ 
&&&&& \downarrow \\ 
\wedge^3 F^* \otimes G^* \otimes S_2(P) & \oplus & 
F^* \otimes S_2(G)^* \otimes P & & & \oplus & 
S_3(P)(-c) & \oplus & 
\wedge^4 F^* \otimes  S_2(P) \\
&&&&& \downarrow \\
\wedge^4 F^* \otimes  S_2(P) & \oplus & 
\wedge^2 F^* \otimes G^* \otimes P & \oplus & 
S_2(G)^*  & \oplus & & & 
F(-c) \otimes S_2(P) \\
&&&&& \downarrow \\ 
&& \wedge^3 F^* \otimes P & \oplus &
F^* \otimes G^* & \oplus & & & 
G(-c) \otimes S_2(P) \\
&&&&& \downarrow \\ 
&&&& \wedge^2 F^*\\
&&&&& \downarrow \\ 
&&&&& I(p) \\
&&&&& \downarrow \\ 
&&&&& 0. 
\end{array} 
$$
With the help of $E_{\bullet}^*$ we get the following beginning for the
resolution of $\ER^4(R/I,R)(-p)$: 
$$
\begin{array}{c@{\ }c@{\ }c@{\ }c@{\ }c@{\ }c@{\ }c@{\ }c@{\ }c@{\ }c@{\
      }c@{\ }c@{\ }c@{\ }c@{\ }c} 
&&& \downarrow \\
\wedge^2 F^*(c)  \otimes P^*  & \oplus & 
G^*(c) \otimes S_2(P)^*  & \oplus & 
F \otimes G \\
&&& \downarrow \\ 
&& F^*(c) \otimes S_2(P)^*  & \oplus &  
S_2(G) \\
&&& \downarrow \\ 
&&& S_3(P)^*(c) \\
&&& \downarrow \\ 
&&& \ER^4(R/I,R)(-p) \\
&&& \downarrow \\ 
&&& 0. 
\end{array} 
$$
Now we apply Proposition \ref{resolution_via_cohomology} and conclude that
in the top row of the resolution of $I(p)$ only the term $S_3 (G)^*
\otimes P$  remains in the minimal resolution because $S_3 (G)^*
\otimes P$ surjects minimally onto $\ER^5(R/I,R)(-p)$. Continuing in this
fashion we  obtain the following resolution:
$$
\begin{array}{c@{\ }c@{\ }c@{\ }c@{\ }c@{\ }c@{\ }c@{\ }c@{\ }c@{\ }c@{\
      }c@{\ }c@{\ }c@{\ }c@{\ }c}
&&& 0 \\ 
&&& \downarrow \\ 
 S_3(G)^* \otimes P  &  & & \oplus & & &  \\ 
&&& \downarrow \\  
F^* \otimes S_2(G)^* \otimes P & & & \oplus & 
S_3(P)(-c) & &  \\
&&& \downarrow \\
\wedge^2 F^* \otimes G^* \otimes P & \oplus & 
S_2(G)^*  & \oplus & & & 
F(-c) \otimes S_2(P) \\
&&& \downarrow \\ 
&& \wedge^3 F^* \otimes P & \oplus &
F^* \otimes G^* & \oplus &  
G(-c) \otimes S_2(P) \\
&&& \downarrow \\ 
&& \wedge^2 F^*\\
&&& \downarrow \\ 
&&& I(p) \\
&&& \downarrow \\ 
&&& 0. 
\end{array} 
$$
 
(ii) Let $t = 2$ and $r=6$. Then $J = J(\psi  \neq I(\psi)$  defines an
arithmetically Cohen-Macaulay scheme. Thus we get as in the previous case
but slightly easier  a resolution
\[
\begin{array}{cccccccccc}
&&& 0 \\
&&& \downarrow \\
S_3 (G)^*  \otimes P & \oplus && S_4(P)(-c_1 ) \\
&&& \downarrow \\
F^* \otimes S_2 (G)^*  \otimes P &&& \oplus &  F(-c_1 ) \otimes S_3 (P) \\
&&& \downarrow \\
\bigwedge^2 F^* \otimes G^*  \otimes P & \oplus & S_2 (G)^*  & \oplus &
G(-c_1 ) \otimes S_3 (P) & \oplus & \bigwedge^2 F (-c_1 ) \otimes S_2 (P) \\
&&& \downarrow \\
\bigwedge^3 F^*  \otimes P & \oplus & F^* \otimes G^*  & \oplus & &&F
\otimes G(-c_1 ) \otimes S_2 (P) \\
&&& \downarrow \\
&& \bigwedge^2 F^*  & \oplus &&& S_2 (G) (-c_1 ) \otimes S_2 (P) \\
&&& \downarrow \\
&&& J(p) \\
&&& \downarrow \\
&&& 0
\end{array}
\]

\end{example}

\begin{remark} \label{discussion_of_minimality} We want to discuss the
  minimality of the resolution 
  described in Theorem \ref{resolution_of_top-dim_part}: \\
(i) By looking at the twists of the free summands occurring in the
resolution above it is clear that for suitable choices of $F, G, P$ and
sufficiently  general maps $\ffi, \psi$ no
further cancellation is possible, i.e., the resolution is minimal. \\
(ii) Let $r$ be even and let $t=1$. 
This case was also studied by Kustin \cite{Kustin}; his main result gives
(up to the degree shifts) the
same resolution as our Theorem \ref{resolution_of_top-dim_part} (cf.\ also
Corollary \ref{one-sect-even-case}).  His
techniques are
completely different from ours, and while they are more complicated, they
in fact
give the maps in the resolution while our techniques use the Horseshoe
Lemma and
hence do not easily give the maps. 

Kustin has proved that his resolution is minimal in the homogeneous case if
the  section does not correspond to a minimal generator of the Buchsbaum-Rim
module.  The latter assumption cannot be removed.  Indeed, the resolution
predicts  that the 
homogeneous ideal of $X = S$ has $r+1$ generators, i.e., that $X$ is an
almost complete intersection because $\codim X = r$. Now consider the
cotangent bundle of $\PP^2$. It has a global section whose zero locus is a
point in $\PP^2$ being a complete intersection. \\
(iii) We suspect that the phenomenon just described is the only instance
that prevents our resolution from being minimal. That means, we hope that
in case  the
resolution of $X$ described above gives the correct number of minimal
generators of $I_X$ then the {\it whole} resolution is minimal. 
\end{remark} 

In the theorem above our focus has been on $X$ rather than on the degeneracy
locus $S$ itself. The interested reader will observe that the same methods
provide a resolution for $S$. 

According to Theorem \ref{summary_for_degeneracy_loci} we know when $X$ is
arithmetically Gorenstein. In this case its minimal free resolution is
self-dual. In order to make this duality transparent we rewrite the
resolution above as follows. 

\begin{corollary} Let $\cBf$ be a Buchsbaum-Rim sheaf of odd rank $r$ and
  first Chern class $c_1$. 
  Let $X$ be the top-dimensional part of  a regular section of
$\cBf$.  
Make the following definitions:
\begin{quote}
\noindent If $i$ is odd, let $\displaystyle \ell = \min \left \{
{\frac{i-1}{2}} ,
{\frac{r-i-2}{2}} \right \}$, and define
\[
A_i = \bigoplus_{j=0}^\ell \bigwedge^{2j+1} F^* \otimes S_{\frac{i-1-2j}{2}}
(G)^*
\]

\noindent If $i$ is even, let $\displaystyle \ell = \min \left \{
{\frac{i}{2}} ,
{\frac{r-i-1}{2}} \right \}$, and define
\[
A_i = \bigoplus_{j=0}^\ell \bigwedge^{2j} F^* \otimes S_{\frac{i-2j}{2}}
(G)^*
\]
\end{quote}
Then $X$ is arithmetically Gorenstein and has a free resolution of the form
$$
0 \rightarrow R(-c_1 )  \rightarrow A_{r-1} \oplus A_1^* (-c_1 ) \rightarrow
A_{r-2} \oplus A_2^* (-c_1 ) \rightarrow \dots \rightarrow
 A_1 \oplus A_{r-1}^* (-c_1 ) \rightarrow I_X \rightarrow 0.
$$
\end{corollary} 

\begin{proof} It follows by Lemma \ref{prop-of-exteriour-powers} that 
$c_1 = c_1(\cBf) = - c_1(\cBf^*)$ is the integer satisfying 
$$
R(-c_1) \cong \wedge^f F^* \otimes \wedge^g G. 
$$
Thus Theorem \ref{resolution_of_top-dim_part} provides the claim. 
\end{proof} 

\begin{example} (i) If the rank of $\cBf$ is 3, this was treated in \cite{Mig-P_gorenstein},
where the following resolution was obtained:
\[
\begin{array}{ccccccccccccccc}
& 0 \\
& \downarrow \\
& R(-c_1 ) \\
& \downarrow \\
G^* & \oplus & F(-c_1 ) \\
& \downarrow \\
F^* & \oplus & G(-c_1 ) \\
& \downarrow \\
& I_X \\
& \downarrow \\
& 0
\end{array}
\]
(ii) If the rank of $\cBf$ is five, then the corollary gives the following
resolution 
\[
\begin{array}{cccccccccccccc}
&&& 0 \\
&&& \downarrow \\
&&& R(-c_1 ) \\
&&& \downarrow \\
 S_2 (G)^* &&   & \oplus & F(-c_1 ) \\
&&& \downarrow \\
F^* \otimes G^* && & \oplus &  G(-c_1 ) & \oplus & \bigwedge^2 F(-c_1 )  \\
&&& \downarrow \\
\bigwedge^2 F^*  & \oplus & G^* & \oplus &&& F \otimes G(-c_1 ) \\
&&& \downarrow \\
&& F^* & \oplus & && S_2 (G)(-c_1 ) \\
&&& \downarrow \\
&&& I_X \\
&&& \downarrow \\ 
&&& 0
\end{array}
\] 
This resolution has been conjectured in \cite{Mig-P_gorenstein}. 
\end{example}

If we consider a regular section of an even rank Buchsbaum-Rim sheaf we
cannot expect to get an arithmetically Gorenstein subscheme as zero
locus. Still the resolution has some symmetry and looks very much like the
corresponding one for the Gorenstein case.  

\begin{corollary} \label{one-sect-even-case} Let $\cBf$ be a Buchsbaum-Rim
  sheaf of even rank $r$ and 
  first Chern class $c_1$. 
  Let $S$ be the zero locus of  a regular section of
$\cBf$.  
Make the following definitions:
\begin{quote}
If $i$ is odd, let $\displaystyle \ell = \min \left \{ {\frac{i-1}{2}} ,
{\frac{r-i-1}{2}} \right \}$, and $\displaystyle \ell' = \min \left \{
\frac{i-1}{2} , \frac{r-i-3}{2} \right \}$.  Define
\[
\begin{array}{c}
A_i = \displaystyle \bigoplus_{j=0}^\ell \bigwedge^{2j+1} F^* \otimes
S_{\frac{i-1-2j}{2}} (G)^* \\ \\
B_i = \displaystyle \bigoplus_{j=0}^{\ell'} \bigwedge^{2j+1} F^* \otimes
S_{\frac{i-1-2j}{2}} (G)^*
\end{array}
\]
\par
\noindent If $i$ is even, let $\displaystyle \ell = \min \left \{
{\frac{i}{2}} ,
{\frac{r-i}{2}} \right \}$, and $\displaystyle \ell' = \min \left \{
\frac{i}{2}
, \frac{r-i-2}{2} \right \}$.  Define
\[
\begin{array}{c}
A_i = \displaystyle \bigoplus_{j=0}^\ell \bigwedge^{2j} F^* \otimes
S_{\frac{i-2j}{2}} (G)^* \\ \\
B_i = \displaystyle \bigoplus_{j=0}^{\ell'} \bigwedge^{2j} F^* \otimes
S_{\frac{i-2j}{2}} (G)^*
\end{array}
\]
\end{quote}
Then $S$ is arithmetically Cohen-Macaulay and has a free resolution of the form
$$
0 \rightarrow R(-c_1 ) \oplus A_{r}  \rightarrow B_1^* (-c_1 ) \oplus
A_{r-1} \rightarrow  B_2^* (-c_1 ) \oplus  A_{r-2} 
 \rightarrow \dots \rightarrow   B_{r-2}^* (-c_1 ) \oplus A_2 \rightarrow A_1  \rightarrow I_S \rightarrow 0.
$$
\end{corollary}


\section{Some applications} 

In the previous section we have seen how much the properties of our
degeneracy loci depend on the properties of the Buchsbaum-Rim sheaf. Now we
will show how this information can be used to construct schemes with
prescribed properties. Moreover, we will explain how sections of the dual
of a generalized null correlation bundle can be studied with the help of our
results. 
\bigskip 

\noindent {\it Construction of arithmetically Gorenstein subschemes
  containing a given scheme} 
\medskip 

Let $X \subset \PP^n$ be an equidimensional  projective subscheme of
codimension $\geq 3$. It 
is rather easy to find 
a complete intersection $Y$ such that $Y$ contains $X$ and both have the
same dimension. The analogous problem where one requires $Y$ to be
arithmetically Gorenstein but not a complete intersection is much more
difficult. This is relevant if one wants to study linkage with respect to
arithmetically Gorenstein subschemes rather than complete intersections. 
We want to explain a solution to this problem. 

Suppose $X$ has codimension $r$. Let us assume that $r$ is odd. Then we
choose a Buchsbaum-Rim sheaf $\Bf$ of rank $r$ on $\PP^n = \Proj R$ given
by an exact sequence 
$$
0 \to \cBf \to \cF \to \cG \to \Mf \to 0. 
$$
For example, we can take the sheafification of the first  syzygy module of an
ideal which is generated by $R$-regular sequence of length $r+1$. 

Next we choose a regular section $s \in H^0(\cBf(j))$ which also belongs to
$H^0_*(\cJ_X \otimes \cF)$. This is possible if $j$ is sufficiently large. 
Let $S$ be the zero-locus of $s$. Then the top-dimensional part $Y$ of $S$
is arithmetically Gorenstein due 
to Theorem \ref{summary_for_degeneracy_loci}. Furthermore,  $s \in H^0_*(\cJ_X
\otimes \cF)$  ensures that $s$ vanishes on $X$. It follows that  $X
\subset Y$ because both are equidimensional schemes of the same dimension.  

For example, it was shown in \cite{Mig-P_gorenstein} that if $\cBf$  is the
cotangent bundle on $\PP^3$, 
twisted by 3, and X is a set of four distinct points, then a section of $\cBf$
can be found vanishing on X and giving a Gorenstein scheme, Y, of degree 5,
whereas the smallest complete intersection containing X has degree 8.

Now assume that $\codim X = r$ is even. Then we take a hypersurface
containing $X$ which is defined by, say $f \in R' = K[x_0,\ldots,x_n]$. Put
$R = R'/f R'$ and choose a Buchsbaum-Rim sheaf $\cBf$ of rank $r-1$ on $Z = \Proj
R$. $X$ has codimension $r-1$ as a subscheme of $Z$. Thus we can find
as in the previous case an arithmetically Gorenstein subscheme $Y \subset
Z$ containing $X$. We can consider $Y$ also as a subscheme of $\PP^n$. As
such it is still  arithmetically Gorenstein, i.e., it has all the
properties we wanted. 
\bigskip

\noindent {\it Construction of $k$-Buchsbaum schemes} 
\medskip 

A projective subscheme $X \subset \PP^n$ is said to be $k$-Buchsbaum for
some non-negative integer $k$ if 
$$ 
(x_0,\ldots,x_n)^k \cdot H^i_*(\cJ_X) = 0 \quad \mbox{for all} \; i \leq \dim
X. 
$$
Note that any equidimensional locally Cohen-Macaulay subscheme is
$k$-Buchsbaum for some $k$. In fact, one should view the notion of a
$k$-Buchsbaum scheme as a refinement of the notion of an equidimensional
locally Cohen-Macaulay subscheme. The idea is to develop for such schemes
a theory  generalizing the  one  for arithmetically
Buchsbaum schemes (cf., for example, \cite{Mig-MM}, \cite{NS}). 

However, so far there are not many examples available where one knows that they
are $k$-Buchsbaum but not $(k-1)$-Buchsbaum. We are going to construct new
examples now. 

Let $R = K[x_0,\ldots,x_n]$ and let $\ffi: R^{n+k}(-1) \to R^k$ be a
homomorphism such that $I(\ffi) = (x_0,\ldots,x_n)^k$. This is true if
$\ffi$ is chosen general enough. A particular choice of such a map is
described in \cite{BV}, p.\ 15. Let us denote the corresponding
Buchsbaum-Rim sheaf by $\cB_k$. Observe that $\cB_1$ is just the cotangent
bundle of $\PP^n$. 

\begin{proposition} \label{k-Buchsbaum_schemes} Let $S$ be the
   degeneracy locus of a morphism  $\psi: \cP
  \to \cB_k$. If $X$ has codimension $n-t+1$ then it holds: If $t = 1$
  or $t=2$ and $n$ is even then $S$ is arithmetically Cohen-Macaulay;
  otherwise  $S$ is $k$-Buchsbaum but not $(k-1)$-Buchsbaum. 
\end{proposition} 

\begin{proof}  The first assertion follows by Theorem
  \ref{summary_for_degeneracy_loci}. 

According to \cite{Buchs-Eisenbud_annihilator} it holds 
$$
\Ann_R S_j(\Mf) = (x_0,\ldots,x_n)^k \quad \mbox{for all} \; j \geq 1. 
$$ 
Hence the second claim is a consequence of Proposition
\ref{coho_of_top-dimensional}. 
\end{proof} 

In case $k = 1$ even more is true. To this end recall that any
arithmetically Buchsbaum subscheme is $1$-Buchsbaum. But the converse is
not true in general. 

Let $N$ be a finitely generated $R$-module. Then the embedding $0 :_N \fm
\hookrightarrow \HH^0(N)$ induces natural homomorphisms of derived functors 
$$
\ffi^i_N: \ER^i(K,N) \to \HH^i(N). 
$$
Due to \cite{SV2}, Theorem I.2.10 $N$ is a Buchsbaum module if and only if the maps
$\ffi^i_N$ are surjective for all $i \neq \dim N$. A subscheme $X \subset
\PP^n$ is called arithmetically Buchsbaum if its homogeneous coordinate
ring $R/I_X$ is Buchsbaum. Now we can show the announced strengthening of
the previous result in case $k=1$. 

\begin{proposition} \label{aBM-schmes} Let $S$ be the
  degeneracy locus of a morphism  $\psi: \cP
  \to \cOP$. If $S$ has codimension $n-t+1$ then it is arithmetically
  Buchsbaum.  
\end{proposition} 

\begin{proof} We will use again the Eagon-Northcott $E_{\bullet}$ complex
  associated to $\psi$. Let $B = H^0_*(\cOP)$ and put $A = R/I_S$. 

We want to show that $\ffi^j_A$ is surjective if $j \neq t = \dim A$. This
is clear if $\HH^j(A)$ vanishes. 
Let $j = n+t- 2i < t$ be an integer such that $\HH^j(A) \neq 0$. 
Using the exact sequences in the proof of Proposition \ref{cohomology-of
  locus} we get diagrams 
$$
\begin{array}{ccc} 
\ER^{n+t-2i}(K,A)(p) & \to & \ER^{n+1-i}(K,\im \delta_i) \\[3pt]
\downarrow {\scriptstyle \ffi^{n+t-2i}_A} & & \downarrow {\scriptstyle
  \ffi^{n+1-i}_{\im \delta_i}}\\[3pt] 
\HH^{n+t-2i}(A)(p) & \to & \HH^{n+1-i}(\im \delta_i) 
\end{array} 
$$ 
and 
$$
\begin{array}{ccc}
\ER^{n+1-i}(K,E_i) & \to & \ER^{n+1-i}(K,\im \delta_i) \\[3pt]
\downarrow {\scriptstyle \ffi^{n+1-i}_{E_i}} & & \downarrow {\scriptstyle
    \ffi^{n+1-i}_{\im \delta_i}}\\[3pt] 
\HH^{n+1-i}(E_i) & \to & \HH^{n+1-i}(\im \delta_i).  
\end{array} 
$$ 
They are commutative because the vertical maps are canonical. Moreover, we
have seen in the proof of Proposition \ref{cohomology-of
  locus} that the lower horizontal maps are isomorphisms. 

Now it is well-known that the modules $\wedge^q B^*$ are Buchsbaum modules
if $1 \leq q \leq n$. Thus the modules $E_q$ are Buchsbaum, too. 
Hence the diagrams show that the surjectivity of
$\ffi^{n+1-i}_{E_i}$ 
 implies this property first for $\ffi^{n+1-i}_{\im \delta_i}$ and then
for $\ffi^{n+t-2i}_A$. It follows  that $A$ is Buchsbaum. 
\end{proof} 

In the special case that $\cP$ has rank $t = n-1$ the last result is also
contained in \cite{Chang-Diff-Geom}. If $t < n-1$ our  result is a little
surprising.   In fact,  the main result of \cite{Chang-Diff-Geom}  has been
generalized in \cite{habil}, Corollary
II.3.3. It says that  arithmetically Buchsbaum subschemes of arbitrary
codimension can be
characterized by means of a particular locally free resolution. As a
consequence, every arithmetically Buchsbaum subscheme of $\PP^n$ is the
zero-locus of a global section of a vector bundle which is the direct sum
of exterior powers of the cotangent bundle.  
\bigskip

\noindent {\it Some vector bundles of low rank and their sections} 
\medskip 

Let $R$ be again a graded Gorenstein $K$-algebra of dimension $n+1$. We
assume that  $n \geq 3$ is an odd integer. 
The aim of this subsection is to show that  a vector
bundle arising from a Buchsbaum-Rim sheaf by quotienting out non-vanishing
sections can be studied by means of our results. Then we construct  vector
bundles of rank $n-1$ on $Z = \Proj R$ and apply this
principle to sections of them. 

Let $\cBf$ be a Buchsbaum-Rim sheaf on $Z$ having  global non-vanishing
sections such there is an exact sequence 
$$
0 \to \cQ \stackrel{\gamma}{\longrightarrow} \cBf \to \cE \to 0 
$$ 
where $H^0_*(Z,\cQ)$ is a free $R$-module of rank $u$ and $\cE$ a vector
bundle on $Z$ of rank $r-u$. 

Now, we want to consider a morphism $\psi: \cP \to \cE$ dropping rank in the
expected codimension $r-u-t+1$. This morphism can be lifted to a morphism
$\beta: \cP \to \cBf$ which provides a morphism $\alpha = (\beta,\gamma): \cP
\oplus \cQ  \to \cBf$. Since the degeneracy locus of  $\gamma$ is empty,
$\cBf^*$ is locally the 
direct sum of $\cE^*$ and $\cQ^*$. It follows that  the images of
$\wedge^t \psi^*$ and $\wedge^{t+u} \alpha^*$  are locally
isomorphic. Hence the degeneracy locus $S$ of $\psi$ and the degeneracy
locus of $\alpha$ agree. Thus $\alpha$ drops  rank in the expected
codimension too and we can apply our previous results. 

Next, we construct  a class of vector bundles which
contains the duals of null correlation bundles. To this end let $I =
(f_0,\ldots,f_n) \subset R$ be a complete intersection. Let $d_i = \deg
f_i$. The first syzygy 
module of $I$ defines a Buchsbaum-Rim module $\Bf$ which fits into the exact
sequence  
$$
0 \to \Bf \to \oplus_{i=0}^n R(-d_i) \stackrel{\ffi}{\longrightarrow} R \to
R/I \to  0. 
$$
The Buchsbaum-Rim sheaf $\cBf = \widetilde{\Bf}$ can often be used to
construct a  vector bundle of rank $n-1$ on $Z$. 

\begin{proposition} \label{new_bundles} Suppose there is an integer $c$
  such that the  degrees satisfy 
$$
c = d_0 + d_1 = d_2 + d_3 = \ldots = d_{n-1} + d_n.  
$$
Then $\cBf(c)$ admits a non-vanishing global section $s$ which gives rise to an
exact sequence  
$$
0 \to \cPP(-c) \stackrel{s}{\longrightarrow} \cBf \to \cN \to 0 
$$ 
where $\cN$ is a vector bundle of rank $n-1$ on $Z$. 
\end{proposition} 

\begin{proof} The Koszul relation of the generators $f_i$ and $f_{i+1}$ of
  $I$ gives rise to a global section of $\cBf(-d_i - d_{i+1})$. Taking the
  sum over these sections with even $i$ yields a section $s$ which does not
  vanish on $Z$ since $I$ is an $\fm$-primary ideal. 
\end{proof} 

In case $Z = \PP^n$ and $d_0 = \ldots = d_n = 1$ the bundle $\cBf$ is 
the cotangent bundle of $\PP^n$ and  $\cN^*$ is
called  {\it
   null correlation bundle} in \cite{OSS} where on  p.\ 79 it is
 constructed in a slightly different way. If $n=3$ then $\cN$ is self-dual.
 Thus we call the dual of a vector bundle constructed as in the proposition
 above {\it generalized null correlation bundle}. Due to  the principle
 described above our previous results apply to multiple sections of the
 dual of a generalized null correlation bundle. We obtain, for example.

\begin{corollary} The degeneracy locus of a multiple section of the dual of
  a  null correlation bundle is arithmetically
  Buchsbaum but not arithmetically Cohen-Macaulay. 
\end{corollary} 

This result is well-known if $n=3$. In fact, in this case the null
correlation bundle can be constructed (via the Serre correspondence) as an
extension  
$$
0 \to \cPP(-2) \to \cN \to \cJ_X \to 0 
$$ 
where $\cJ_X$ denotes the ideal sheaf of two skew lines in $\PP^3$ (cf.\
\cite{Barth-1977}, p. 145 or \cite{Ellia-Fiorentini}).


\begin{thebibliography}{999} 

\bibitem{Barth-1977} W.\ Barth, {\em Some properties of stable rank-$2$
    vector bundles on $\PP^n$}, Math.\ Ann. {\bf 226} (1977), 125--150. 

\bibitem{Bruns-Herzog} W.\ Bruns, J.\ Herzog, {\em ``Cohen-Macaulay
    rings''}, Cambridge Studies in Advanced Mathematics {\bf 39}, Cambridge
  University Press, 1993.

\bibitem{BV} W.\ Bruns, U.\ Vetter, {\em ``Determinantal rings''}, LNM
  1327, Springer-Verlag, 1988. 

\bibitem{Buchs-Eisenbud_annihilator} D.\ A.\ Buchsbaum, D.\ Eisenbud, {\em
    What annihilates a module?}, J.\ Algebra {\bf 47} (1977), 231--243. 

\bibitem{Chang-Diff-Geom} M.\ Chang, {\em Characterization of
    arithmetically Buchsbaum subschemes of codimension $2$ in $\PP^n$},
  J.\ Differential Geom. {\bf 31} (1990), 323--341. 

\bibitem{Eisenbud-Buch} D.\ Eisenbud, {\em ``Commutative algebra with a
    view toward algebraic geometry''}, Graduate Texts in Math. {\bf 150},
  Springer-Verlag, 1995. 

\bibitem{Ellia-Fiorentini} Ph.\ Ellia, M. Fiorentini, {\rm Quelques
  remarques sur les courbes arithmetiquement  Buchsbaum de l'espace
  projectif},  Ann.\ Univ.\ Ferrara {\bf XXXIII} (1987), 89-111.

\bibitem{KMNP} M.\ Kreuzer, J.\ C.\ Migliore, U.\ Nagel, C.\ Peterson,
{\em Determinantal schemes and Buchsbaum-Rim sheaves}, Preprint, 1997.

\bibitem{Kustin} A.\ Kustin, {\em The minimal free resolution of the
    Migliore-Peterson Rings in the case that the reflexive sheaf has even
    rank}, Preprint, 1996. 

\bibitem{Mig-MM} J.\ Migliore, R.\ M.\ Miro-Roig, {\em On $k$-Buchsbaum
    curves}, Comm.\ Algebra {\bf 18} (8) (1990), 2403--2422. 

\bibitem{Mig-P_gorenstein} J.\ Migliore, C.\ Peterson, {\em A
    construction of codimension three arithmetically Gorenstein subschemes
    of projective space}, Trans.\ Amer.\ Math.\ Soc.\ (to appear). 

\bibitem{NS} U.\ Nagel, P.\ Schenzel, {\em Cohomological 
annihilators and Castelnuovo-Mumford regularity}, Contemp.\ Math. {\bf 159}
(1994), 307--328. 

\bibitem{habil} U.\ Nagel, {\em On arithmetically Buchsbaum subschemes and
    liaison}, Habilitationsschrift, 1996. 

\bibitem{OSS} C.\ Okonek, M.\ Schneider, H.\ Spindler, {\em ``Vector
    bundles on complex projective spaces''}, Birkh\"auser, 1980. 

\bibitem{R1} P.\ Rao, {\it Liaison among Curves in} $\PP^{3}$, Invent.\ Math.\
{\bf 50} (1979), 205--217.

\bibitem{S} P.\ Schenzel, {\em ``Dualisierende Komplexe in der lokalen
    Algebra und Buchsbaum-Ringe''}, LNM {\bf 907}, Springer-Verlag, 1982.

\bibitem{SV2}  J.\ St{\"u}ckrad, W.\ Vogel, {\em ``Buchsbaum rings and
    applications''}, Springer-Verlag, 1986.

\bibitem{weibel} C.\ Weibel, {\em ``An introduction to homological algebra''},
Cambridge University Press, Cambridge studies in advanced mathematics 38
(1994).



\end{thebibliography}
\end{document}